\documentclass[10pt]{article}
\usepackage{geometry}
\usepackage{cancel,amssymb,amsmath,mathrsfs,upgreek,mathtools,subcaption,textcomp,ulem,mathdots,color,tikz,tikz-cd,todonotes,framed}
\usetikzlibrary{patterns, decorations.markings,arrows, calc}
\usepackage[colorlinks=true, pdfstartview=FitV, linkcolor=darkblue, citecolor=darkblue, urlcolor=darkblue]{hyperref}
\definecolor{shadecolor}{rgb}{0.9, 0.9, 0.86}
\definecolor{darkgreen}{rgb}{0.2, 0.5,  0}
\definecolor{darkblue}{rgb}{0.1,0.1,0.45}
\definecolor{red}{rgb}{0.9,0,0}

\def\XXint#1#2#3{{\setbox0=\hbox{$#1{#2#3}{\int}$ }
\vcenter{\hbox{$#2#3$ }}\kern-.6\wd0}}

\def\thetadiv{{(\Theta)}}
\def\ss{{\mathbf{s}}}

\def \Cartan#1{\vartheta{#1}}
\def \bea#1\eea {\begin{align} #1 \end{align}}
\def \moduli{ \mathcal{V}}
\def \charvar{\mathfrak R(\mathcal C, \GL_n)}

\def \Berg{\mathbf{B}}
\def\&{\hspace{-15pt}&}
\def\nn{\nonumber}
\def \br{\begin{remark}\rm }
\def \er{\hfill $\triangle$\end{remark}}
\def\liouville{\lambda}
\def \scr{\mathscr}
\def \Malg{\Xi}
\def \u{\mathfrak A}
\def\nf{P}
\def \bs {\boldsymbol}
\def \K{ \mathscr K}
\def \P {\mathbb P}
\def \conn{\scr A}

\def \eqref#1{(\ref{#1})}

\def \wt{\widetilde}
\def\Cauchy{ \mathbf {C}}
\def \wh{\widehat}
\newcommand{\C}{\mathbb{C}}
\renewcommand{\le}{\left}
\newcommand{\ri}{\right}

\newcommand{\1}{\mathbf{1}}

\newcommand{\g}{\gamma}
\newcommand{\E}{\mathrm{E}}
\renewcommand{\b}{\beta}
\renewcommand{\a}{\alpha}

\renewcommand{\d}{\mathrm d}
\renewcommand{\O}{\mathcal{O}}
\renewcommand{\S}{\Sigma}

\let\ii\i
\renewcommand{\i}{\mathrm{i}}
\newcommand{\pa}{\partial}
\def\res{\mathop{\mathrm {res}}\limits_}
\def \tr {\mathrm{tr}\,}
\def\be{\begin{equation}}
\def\ee{\end{equation}}

\def\SL{\mathrm{SL}}
\def\GL{\mathrm{GL}}
\def\Mat{\mathrm{Mat}}

\newtheorem{theorem}{Theorem}[section]
\newtheorem{theoremdefinition}[theorem]{Theorem-Definition}
\newtheorem{example}[theorem]{Example}

\newtheorem{lemma}[theorem]{Lemma}
\newtheorem{remark}[theorem]{Remark}

\newtheorem{proposition}[theorem]{Proposition} 
\newtheorem{corollary}[theorem]{Corollary}
\newtheorem{definition}[theorem]{Definition}

\def\framing{{\mathfrak{S}}}

\def\QED {\hfill $\blacksquare$\par\vskip 10pt}

\def\CC{\mathcal{C}}
\def\E{\mathcal{E}}
\def \DD{\mathbb D}
\def\F{\mathbf{F}}
\def\g{\gamma}

\renewcommand{\theequation}{\arabic{section}.\arabic{equation}}

\makeatletter
\@addtoreset{equation}{section}
\makeatother

\begin{document}

\vspace{0.2cm}
\begin{center}
\begin{Large}
\textbf{
Higgs fields, non-abelian Cauchy kernels and the Goldman symplectic structure} 
\end{Large}
\end{center}

\begin{center}
M. Bertola$^{\dagger\ddagger\bullet\star}$ \footnote{Marco.Bertola@\{concordia.ca, sissa.it\}} 
C. Norton$^{\clubsuit}$ \footnote{nchaya@umich.edu}
G. Ruzza $^{\diamondsuit\circ}$ \footnote{gruzza@fc.ul.pt}
\\
\bigskip
\begin{minipage}{0.7\textwidth}
\begin{small}
\begin{enumerate}
\item [${\dagger}$] {\it  Department of Mathematics and Statistics, Concordia University\\ 1455 de Maisonneuve W., Montr\'eal, Qu\'ebec, Canada H3G 1M8} 
\item[${\ddagger}$] {\it SISSA, International School for Advanced Studies, via Bonomea 265, Trieste, Italy }
\item[$\bullet$] {\it Istituto Nazionale di Fisica Nucleare (INFN), sezione di Trieste}
\item[${\star}$] {\it Centre de recherches math\'ematiques, Universit\'e de Montr\'eal\\ C.~P.~6128, succ. centre ville, Montr\'eal, Qu\'ebec, Canada H3C 3J7}
\item[${\clubsuit}$] {\it Department of Mathematics, University of Michigan, Ann Arbor, MI 48109, USA}
\item[${\diamondsuit}$] {\it Grupo de F\'\ii sica Matem\'atica, Departamento de Matem\'atica, Instituto Superior T\'ecnico, Av. Rovisco Pais, 1049-001 Lisboa, Portugal}
\item[${\circ}$] {\it Departamento de Matem\'atica, Faculdade de Ci\^encias da Universidade de Lisboa, Campo Grande Edif\'\ii cio C6, 1749-016, Lisboa, Portugal}
\end{enumerate}
\end{small}
\end{minipage}
\end{center}

\vskip 15pt

\begin{center}
\begin{abstract}
We consider the moduli space of vector bundles of rank $n$ and degree $ng$ over a fixed Riemann surface of genus $g\geq 2$ with the explicit parametrization in terms of the Tyurin data.
The  ``non-abelian'' theta divisor consists of bundles such that~$h^1\geq 1$. 
On the complement of this divisor we construct a non-abelian (i.e. matrix) Cauchy kernel explicitly in terms of the Tyurin data. With the additional datum of a non-special divisor, we can construct a reference flat holomorphic connection which
also depends holomorphically on the moduli of the bundle. This allows us to identify the bundle of Higgs fields, i.e. the cotangent bundle of the moduli space, with the affine bundle of holomorphic connections and provide a monodromy map into the $\GL_n$~character variety. We show that the Goldman symplectic structure on the character variety pulls back along this map to the complex canonical symplectic structure on the cotangent bundle and hence also on the space of  connections.
The pull-back of the Liouville one-form to the affine bundle of connections is then shown to be a  logarithmic form with poles along the non-abelian theta divisor and residue given by~$h^1$. 
\end{abstract}
\end{center}

\section{Introduction and results}

In the study of stable holomorphic vector bundles over a Riemann surface $\CC$ of genus $g>0$ the result of principal importance is the theorem by Narasimhan and Seshadri \cite{NaraSesh}. It states that a holomorphic vector bundle  is stable if and only if it comes from an irreducible projective unitary representation  of the fundamental group of the surface. Donaldson showed that the condition is equivalent
to the existence of a unitary connection with curvature taking values in the center of the group \cite{Donaldson}. The connection can be obtained from the  minimization of a functional and depends in a  smooth way on the moduli of the bundle, but not  {\it holomorphically}. 
This famous result establishes a one-to-one correspondence between stable holomorphic bundles and 
projective unitary representations and thus provides the moduli space with the structure of real smooth manifold; trace invariants of monodromies are the smooth functions on the moduli space. 

Of  a similar nature is  the even more classical result of Fuchsian uniformization of Riemann surfaces, which provides real coordinates on the moduli space of Riemann surfaces of genus $g\geq 2$; the Schwarzian derivative with the uniformizing coordinate provides a natural {\it projective connection} on the Riemann surface $\CC$. However, again, this projective connection, while holomorphic with respect to the conformal structure of the Riemann surface,  does not depend holomorphically on the moduli.  

The space of (holomorphic) projective connections is an affine space modelled on the space of quadratic differentials on $\CC$; the latter is identified with the  fiber of the cotangent bundle  to  the Teichm\"uller space $\mathcal T_g$ above $\CC$. 
A projective connection is equivalent to the datum of an ``oper'', namely, a second order ordinary differential equation 
\be
\label{proj0}
(\pa_z^2 + \scr U(z))\psi(z)=0
\ee
for a holomorphic $(-1/2)$differential $\psi$, where $\scr U$ is the projective connection.  The monodromy representation of \eqref{proj0} is a  function of the chosen projective connection $\scr U$. 
The set of holomorphic projective connections is an affine bundle over $\mathcal T_g$. The fiber is an affine space modelled over the corresponding fiber of the cotangent bundle and, hence, over the space of holomorphic quadratic differentials. 
If we choose a reference projective connection $\scr U_\CC$ then we can identify the cotangent space with the space of projective connections by the map $\phi \mapsto  \scr U_\CC + \phi$. If the reference connection $\scr U_\CC$ depends holomorphically on the moduli, then the above map allows us to map 
 $T^* \mathcal T_g$  into the affine bundle of projective connections and then to the $\P\SL_2$  character variety  by the monodromy map.

On the $\P\SL_2(\C)$ character variety there is a natural complex-analytic Poisson symplectic structure due to Goldman \cite{Goldman}; the question posed and answered in \cite{BKN,Kawai, Tak} was precisely  how to choose $\scr U_\CC$ in such a way that the above map from $T^*\mathcal T_g$ to the character variety is complex-analytic and establishes a  Poisson morphism between the canonical complex symplectic structure on $T^*\mathcal T_g$ and the character variety. For this purpose,  the Fuchsian projective connection is not appropriate because it does not depend analytically on the moduli. 
\vskip 10pt
A parallel problem can be posed in the context of holomorphic vector bundles. 
The space of holomorphic (flat) connections on a flat holomorphic vector bundle $\mathcal X$  of rank $n$ is also an affine vector space modelled on the space of Higgs fields $ H^0({\rm End}( {\mathcal X}) \otimes \mathcal K)$, where $\mathcal K$ denotes the canonical bundle of the curve $\CC$. 
Similarly to the case  of quadratic differentials, a holomorphic Higgs field on a bundle $\mathcal X$ can be thought of as a cotangent vector of the  holomorphic cotangent bundle of the moduli space of vector bundles. 
Hence we can view the space, $\scr A$, of holomorphic connections  as an affinization of the cotangent bundle of the moduli space. 
The main thrust in this direction came in \cite{krichever, kricheveriso}; in \cite{krichever}  
the author describes an open-dense locus in the space of Higgs fields in terms of the Tyurin parametrization  whereby the Higgs field can be viewed as a matrix-valued meromorphic differential with appropriate conditions on the structure of singularities. 
This allows one to write the canonical symplectic two-form on the cotangent bundle of the moduli space as a residue pairing localized at the points of the Tyurin divisor.
In the subsequent \cite{kricheveriso} the author takes a similar approach by considering meromorphic connections described in the same framework. He then takes the same formula for the symplectic structure and translates it in terms of the entries of the (generalized) monodromy matrices. The result is an explicit two-form on the space of monodromy data. We are going to comment more on this after Theorem \ref{thmmalgfay} in Remark \ref{remkri}.
\\[1pt]

With regard to symplectic structures, the broad strokes picture is then quite similar to the situation with opers and projective connections;  on one side  the $\P \GL_n(\C)$ character variety  has a complex-analytic Poisson symplectic structure due to Goldman, and on the other side the holomorphic cotangent bundle of the moduli space has a  natural complex symplectic structure. 
In the context of opers and projective connections it was shown in  \cite{BKN}  that  the monodromy map translates to a Poisson morphism {\it only for certain classes} of projective connections that depend holomorphically on the moduli of the curve; these include the so-called Bers, Schottky, Bergman and Wirtinger projective connection and exclude the Fuchsian one since it is not holomorphic in the moduli.

The question we consider here is the exact analog of the case of opers  considered in \cite{BKN, Kawai, Tak} and extending the considerations of  \cite{kricheveriso} in the case of vector bundles; namely, whether there is a natural choice of reference holomorphic connection for every vector bundle $\mathcal X$, depending holomorphically on the moduli of the bundle, and such that the mapping of the cotangent bundle to the character variety results in a Poisson morphism on the complement of the non-abelian theta divisor.  
In \cite{kricheveriso} the Poisson morphism was considered from the space of connections itself\footnote{
It should be mentioned that the setting of loc. cit. includes also the generalized monodromy data (i.e. Stokes' matrices) for connections of arbitrary Poincar\'e~rank.} (which was shown to be an affinization of the space of Higgs fields) but without an explicit choice of reference connection; the relationship between the bundle of Higgs fields \cite{krichever} and the bundle of connections \cite{kricheveriso} requires to explicitly fix an origin in the affine space. Such a choice can be considered as the bridge between the two works \cite{krichever, kricheveriso}. Only after such reference point is chosen it  is then meaningful to consider the morphism as a map from the cotangent bundle into the character variety. 

In this work we introduce  a global meromorphic section with first order poles on a divisor as a reference connection. 

Our central object of interest is the Liouville one form rather than the symplectic form; we show that the  singularity of its pull-back along the theta divisor is a  simple pole with positive integer residue and that this integer has a clear cohomological interpretation  (see Theorem \ref{theoremFay} and below for more details and illustrative example). This  can be viewed as a statement about the integrality of the symplectic form.

\paragraph{Outline of results.}
We consider the moduli space $\moduli$  of  vector bundles $\E$ of rank $n$ and  degree $ng$ such that the generic fiber is spanned by the holomorphic sections; these vector bundles admit a convenient and explicit parametrization due to Tyurin \cite{Tyurin1,Tyurin2}.
The notion of Tyurin data was notably used in \cite{EnRub} and, as mentioned,  \cite{krichever, kricheveriso}.

We first re-derive a self-contained  treatment of Tyurin data in Section~\ref{Tyurin} including the general case of arbitrary multiplicity of the Tyurin points. We work mostly in the context of {\it framed} vector bundles; namely we choose a point $\infty\in \CC$ and a basis of global holomorphic sections that trivialize the fiber at $\infty$.

In terms of the Tyurin data we can derive the Riemann--Roch theorem for vector bundles in the same spirit as the classical treatment of line bundles, using an analogue of the Brill--Noether matrix, see Section \ref{BrillNoether}. \par \vskip 4pt

The central object in our construction is the non-abelian Cauchy kernel for a framed vector bundle $\E$, a generalization of the classical notion of scalar Cauchy kernel on a Riemann surface \cite{BehnkeStein,Fay73,GusmanRodin,Rodin,Zvero}.
This kernel is a matrix, $\Cauchy_\infty({q},{p})$ whose rows (as functions of $p$) are meromorphic sections of the bundle $\E$ with pole at $q$ and zero at $\infty$ while the columns (as functions of ${q}$) are meromorphic sections of $\E^*\otimes \mathcal K$ with only poles at $p,\infty$.
By using a Riemann--Hilbert factorization approach to vector bundles on a Riemann surface~\cite{Rodin}, we interpret sections of these bundles as vector-valued functions or differentials with suitable conditions on their poles and zeros, determined from the Tyurin data.
The non-abelian Cauchy kernel is therefore uniquely determined as a matrix-valued function/differential on the Riemann surface by controlling its poles relative to the Tyurin data and by the normalization condition that $\res{q=p} \Cauchy_\infty(q,p)$ is the identity matrix, see Section \ref{secCauchy} and Theorem-Definition \ref{defCauchy}.

Similarly to the scalar case, this non-abelian Cauchy kernel exists  on the locus $h^1(\E)=0$ and has a simple pole as a function of the holomorphic moduli on the non-abelian theta divisor, i.e., the locus $\thetadiv\subset \moduli$ where $h^1(\E)>0$.
We show how to construct {\it explicitly} such kernels in term of the Tyurin data, see \eqref{matrix}  in the general case.
For the reader with some familiarity with the Tyurin parametrization, at least in the generic case, we offer here the following explicit description.
Suppose all Tyurin points $t_1,\dots,t_{ng}$ in the Tyurin divisor $\scr T$ are simple, and  let $\mathbf {h}_j$ be the Tyurin vectors as in \eqref{TyurinVectors} (these vectors were denoted by $\alpha$ in \cite{krichever}).Then there exists a unique $n\times n$ {\it (matrix) Cauchy kernel} $\Cauchy_\infty({q},{p})$ such that 
\begin{itemize}
\item[-] it is a meromorphic matrix--valued  differential with respect to the variable ${q}$ with its divisor satisfying  $(\Cauchy_\infty({q},{p}))_{q}\geq -{p}-\infty$;
\item[-] it is a matrix--valued meromorphic function with respect to  the variable ${p}$ with its divisor satisfying $(\Cauchy_\infty({q},{p}))_{p}\geq \infty-{q}-\scr{T}$,
\end{itemize}
and such that $\mathbf h_j \Cauchy_\infty(p_j,p)=0,\ \forall p_j\in \scr T$, together with the normalization that the residue along the diagonal is the identity.

We give an explicit expression of the Cauchy kernel in terms of the Tyurin data as ratio of determinants:  in \eqref{expressioncauchyblock} for the case of arbitrary Tyurin data  and Example \ref{exTyusymp} for the generic case of simple Tyurin points. These formulas generalize similar classical formulas for scalar Cauchy kernels.

The  Cauchy kernel has the following two primary uses:
\begin{enumerate}
\item[-] it allows us to construct all Higgs fields explicitly (see Section~\ref{secHiggs}); 
{ for example the expression
\be
\label{omegaphiintro}
\Phi({q}):= \sum_{t\in\scr T}\res{{p}=t}\Cauchy_\infty({q},{p}) \nf(p) \phi_{t}(p)\nf^{-1}(p)
\ee
gives an arbitrary (framed) Higgs field in terms of local holomorphic data at the Tyurin divisor. Here $P$ is the polynomial normal form of the transition function of the bundle as defined in Section~\ref{Tyurin}. 
 }  
\item[-] {the regular part along the diagonal } defines a connection that depends holomorphically on the moduli of the bundle (Section~\ref{secAffine}). This allows to holomorphically map the cotangent bundle of the moduli space to the character variety.
To construct this connection we fix an arbitrary reference line bundle of degree $g$ corresponding to a non-special  divisor $\scr D$ and tensor $\E$ by its dual so as to obtain a bundle of degree zero.
Then we define the reference holomorphic connection ${\d-}\F_\scr D$ on this resulting bundle by \eqref{defFDelta}. 
\end{enumerate}

\begin{figure}
\begin{equation*}
\begin{tikzcd}[ampersand replacement=\&]
 T^*\moduli_0  
\arrow{dr} 
\arrow[r,"{\mathbf A}= \F_\scr D + \Phi"]
\&
\conn^{\scr D}   \arrow [d]  \arrow [r,"\mathcal M"] \& \charvar
\\[20pt]
 \& \moduli_0 \arrow [swap, from=2-2, to=1-2, bend right, "\F_\scr D", hookrightarrow, right]
\end{tikzcd} 
\end{equation*}
\caption{Illustration of the main spaces and maps.}
\label{figintro}
\end{figure}

The connection $\d -\F_{\scr D}$ is holomorphic on the bundle $\E^*\otimes\mathcal O(\scr D)$ or can be equivalently viewed as meromorphic on the trivial vector bundle {following ideas first presented in \cite{kricheveriso}}; in this latter perspective $\F_{\scr D}$ is simply a matrix of meromorphic differentials on $\CC$  such that the local monodromy  of the connection around each pole is trivial (Theorem~\ref{nosing}). This induces a map, parametrized by the choice of $\scr D$, from the moduli space $\moduli$ to the $\GL_n(\C)$ character variety of the fundamental group  of $\CC$ 
\be
\charvar:={\rm Hom}\big( \pi_1(\CC, \infty) , \GL_n\big)/\P\GL_n,
\ee
where the action of $\P\GL_n$ is the global conjugation.

The choice of the divisor  $\scr D$ affects only projectively the representation  (Section~\ref{secPGLn})  which is therefore independent of $\scr D$  in  the  $\P \GL_n$ character variety: this means that we can construct a canonically defined immersion of the moduli space $\moduli$  in the $\P \GL_n$ character variety; since $\Cauchy_\infty$ depends holomorphically on the Tyurin data, so does the representation. This construction provides an   alternative description  to the projective-unitary one used in conjunction with Narasimhan--Seshadri's theorem.\par\vskip 10pt

In Section~\ref{secAffine} we then consider  $\conn^{\scr D}$,  the bundle  over $\moduli$ whose fiber consists of connections that are holomorphic on $\E^*\otimes \mathcal O(\scr D)$: the reference connection  $\F_{\scr D}$  provides a section of the canonical projection $\conn^{\scr D} \to \moduli$  which is holomorphic away from the non-abelian theta divisor $\thetadiv$ and with a simple pole therein. This section allows us to identify the cotangent bundle $T^*\moduli$ (away from the theta divisor) with $\conn^{\scr D}$ by sending the Higgs field $\Phi$ to the connection $\d - (\F_{\scr D} + \Phi)$. Composing  this identification with the monodromy map, we obtain a map from the cotangent bundle $T^*\moduli$ to the character variety $\charvar$.
The resulting picture is illustrated in Fig. \ref{figintro}.\par \vskip 10pt

We examine the symplectic properties of the monodromy map in Section~\ref{secEquiv}. 
Consider  the  tautological holomorphic (Liouville) one-form $\lambda$ on the (holomorphic) cotangent bundle of the moduli space:  we can pull it back onto the space of connections $\conn^{\scr D}$ by our choice of reference connection $\F_{\scr D}$.
The resulting expression $\Malg$ is the {\it Malgrange--Fay one-form} (Definition \ref{defMalFay}), written as an integral over the canonical dissection of the curve $\CC$. This defines a symplectic structure on $\conn^{\scr D}$. 
On the other hand, we can pull back to $\conn^{\scr D}$ the (complex) symplectic structure induced by the 
Goldman Poisson bracket \cite{Goldman}  (which is non-degenerate for closed surfaces). 
Then our result is that these two structures coincide up to a simple proportionality factor (Theorem~\ref{main}).

Finally we analyze the singularity structure of the Malgrange--Fay form $\Malg$; it is ill-defined on the non-abelian theta divisor $\thetadiv$ and we show that it  is a {\it logarithmic form}  with a simple pole along $\thetadiv$ and  residue equal to $ h^1(\E)$ (Theorem~\ref{theoremFay}):
$$
{
\res{\E\in \thetadiv} \Malg = h^1(\E).}
$$
 This is the parallel of a  result of Fay (Theorem 2 in \cite{Fay92}). However the form that Fay introduced in loc. cit.  was not identified with the tautological form on the cotangent bundle of the moduli space, and furthermore no connection with the Goldman bracket was established. 

As a logarithmic form, the Malgrange--Fay form $\Malg$ has the same singularities, locally, as the differential of the logarithm of a determinant; this is identified as the determinant of the inverse to the Brill--Noether--Tyurin matrix used in the construction of the non-abelian Cauchy kernel. 

\paragraph{Example: line bundles.}
This example illustrates the gist of the results of our paper and it will be generalized to the vector bundle case.

Given a generic line bundle $\E$ (equivalently, a non-special divisor $\scr T$) of degree $g$ and a point $\infty\in \CC$ there is a unique Cauchy kernel $\Cauchy_\infty(q,p)$ with poles at $\scr T, q$ and a zero at $\infty$ in the variable $p$ and poles at $\scr T, p,\infty$ in the variable $q$. It is a differential with respect to  $q$ and a function with respect to $p$. The normalization is determined by $\res{q=p}\Cauchy_\infty(q,p)=1$.

In this case the reference connection $\d-\F_{\scr D}$ can be written simply in terms of Riemann theta functions as follows:
\be
\F_{\scr D} =\d_p\ln \frac {\Theta(\u(p-\infty)-{\bf f})}{\Theta(\u(p-\scr D)-\K)} +\vec \omega(p)\cdot  \nabla\ln \Theta({\bf f}),
\ee
where $\u$ denotes the Abel map, $\K$ the vector of Riemann constants, 
 ${\bf f} = \u(\scr T-\infty)+\K$, $\vec \omega$ is the row vector of the  $g$ holomorphic $\a$-normalized differentials and $\nabla$ denotes the gradient operator.
Inspection shows that $\F_{\scr D}$ has a simple pole at every point ${\bf f}$ of the Jacobian where $\Theta=0$. These correspond precisely to line bundles  $\scr T$ of degree $g$ such that the corresponding line bundle $\E$ has $h^1(\E)\geq 1$. The connection $\F_{\scr D}$ is  a holomorphic connection  on $\mathcal O(\scr T-\scr D)$ and is single-valued on the Picard variety (i.e. with respect to $\scr T$ or, equivalently, ${\bf f}$ in the Jacobian). 

The fiber of $\conn^{\scr D}$ above $\E$  can be identified with meromorphic differentials of the form ${\bf A}= \phi + \F_\scr D$ where $\phi$ (a holomorphic differential) represents the Higgs field. 
In this case the character variety consists of the exponentials of the  periods of ${\bf A}$ (note that the integer residues at the points of $\scr T, \scr D$ do not contribute to the monodromy representation). This gives a (scalar) representation of the fundamental group with monodromy factors
\be
\chi(\gamma) = \exp\oint_{\gamma} {\bf A}.\nn
\ee
A direct computation of the Malgrange form (Def. \ref{defMalFay}) in this simple case gives $\Xi = -{\bf c}^t \cdot \delta {\bf f}$ where ${\bf c}$ is the vector of coefficients of $\phi = \sum c_\ell \omega_\ell$.
Here and in the rest of the paper $\delta$ denotes the exterior derivative operator on $\moduli$ In this scalar case the Goldman Poisson bracket is simply
\be
\nn
\{\ln \chi(\gamma), \ln \chi(\wt \gamma)\} = \gamma \circ \wt \gamma,
\ee
denoting by $\circ$ the intersection number of the curves.
This is a non-degenerate Poisson bracket and corresponds to the symplectic form 
\be
\sum_{j=1}^g\delta \ln \chi (\alpha_j) \wedge \delta \ln \chi (\beta_j)=\sum_{j=1}^g\delta  \oint_{\alpha_j} {\bf A} \wedge \delta \oint_{\beta_j} {\bf A} =  -2\pi\i\, \delta {\bf f} \wedge \delta {\bf c} = -2\pi\i\, \delta \Malg.
\ee
In order to pull back the corresponding symplectic form on $\conn^{\scr D}$ and see the type of singularity along the theta divisor we need to choose a locally smooth section   of the map $\conn^{\scr D}\to \moduli_0$: we choose
${\bf A} =\F_{\scr D} + \vec \omega(p)({\bf c}- \nabla \ln \Theta({\bf f}) )$. Then  the Malgrange form becomes 
$$\Xi= \delta \ln \Theta({\bf f}) - {\bf c}^t \cdot \delta {\bf f}.
$$
This formula shows its logarithmic nature on the theta divisor.
The fact that the residue of $\Malg$ at a point $\E$ of  $\thetadiv$ is $h^1(\E)$ is simply a restatement of the classical Riemann's  singularity  theorem.

\paragraph{Common notations}
\begin{itemize}
\item $\CC$: a curve of genus $g\geq 2$, with a choice of point $\infty\in\CC$. $\mathcal K$: the canonical bundle of $\CC$
\item $\moduli$: the moduli spaces of vector bundles, $\E$, of rank $n$ and  degree $ng$  such that the fiber over $\infty$ is spanned by global holomorphic sections (as is the generic fiber).
\item $\delta$: exterior differential on~$\moduli$.
\item $\wh \moduli$: like the above but also with $n$ ``framing'' sections linearly independent at $\infty$.
\item  $\thetadiv, \wh \thetadiv$: the non-abelian theta divisors where $h^1(\E)>0$ in $\moduli, \wh \moduli$, respectively.
\item $\moduli_0:= \moduli \setminus \thetadiv$, $\wh \moduli_0:= \wh \moduli \setminus \wh \thetadiv$.
\item $\scr T$ the Tyurin divisor of a vector bundle see Section~\ref{Tyurin}. We use $T$ for the Tyurin space \eqref{tyurinspace}.
\item $\mathbb T$: the Brill--Noether--Tyurin matrix \eqref{matrix}. The theta divisor is the locus $\det \mathbb T=0$.
\item $\conn^{\scr D}_\E$: the affine space of holomorphic connections on $\E^*\otimes \mathcal O(\scr D)$. Then $\conn^{\scr D}$ is the corresponding bundle over $\moduli$, whose fiber at $\E$ is $\conn^{\scr D}_\E$.  See Section~\ref{secAffine}.   
\item $\omega_{j},\ j=1,\dots, g$ a Torelli marked basis of holomorphic differentials; $\omega_{z,w}(p)$ the third kind normalized differential with residue $+1$ at $z$ and $-1$ at $w$. 
\end{itemize}

\paragraph{Note added.}
After completion of this work, a parallel result was proved in~\cite{TakArx} using local holomorphic sections to induce a global holomorphic symplectic form on the bundle of connections.
In both cases, the monodromy map is a symplectomorphism for the Goldman bracket and therefore the symplectic structures are {\it equivalent} in the sense (see~\cite{BKN}) that they induce the same (complex) symplectic structure on the space of connections.
This raises the question of what is the generating function for the corresponding change of polarization: in the context~\cite{BKN} of projective connections such a generating function was conjectured in loc. cit. to be the Selberg zeta-function corresponding to the quasi-Fuchsian group, a conjecture recently proved in~\cite{ErikssonEtAl}.

\paragraph{Acknowledgements.}
The authors wish to thank Prof. I. Krichever for detailed explanations of his prior works on the subject and suggestions for improvement. We also thank  Prof. I. Biswas for discussion related to stability. 
This material is based upon work supported by the National Science Foundation under Grant No. 1440140, while the first author was in residence at the Mathematical Sciences Research Institute in Berkeley, California, during the fall semester of 2019.

C.N.  was supported from the AMS-Simons Travel Grants, which are administered by the American Mathematical Society with support from the Simons Foundation.

G.R. acknowledges support from European Union's H2020 research and innovation programme under the Marie Sk\l owdoska--Curie grant No. 778010 {\it IPaDEGAN}, from the Fonds de la Recherche Scientifique-FNRS under EOS project O013018F, and from the FCT grant 2022.07810.CEECIND.

\section{Tyurin parametrization of vector bundles}
\label{Tyurin}

Throughout this paper, let $\CC$ be a fixed compact Riemann surface of genus $g\geq 2$.
For any vector bundle $\E$ of rank $n$ and degree $ng$ on $\CC$, the Riemann--Roch theorem
\be
\label{riemannroch}
h^0(\E)= h^1(\E) + \deg \E-n(g-1) =  h^1(\E) +n, \qquad h^i(\mathcal{E}):=\dim H^i(\E),
\ee
implies that $\E$ has at least $n$ holomorphic sections. We shall see in Section \ref{BrillNoether} that generically, we have $h^0(\E)=n$, namely $\E$ has exactly $n$ holomorphic sections.
The degree $ng$ case is particularly interesting, as the $n$ sections provide a convenient parametrization of (an open part in) the moduli space of vector bundles on $\CC$, as we now explain expanding on the works of Tyurin and Krichever \cite{kricheveriso,Tyurin1,Tyurin2}.

\subsection{Moduli of framed vector bundles}
Let us fix a point $\infty\in\CC$, and let $\moduli$ be the set of isomorphism classes of holomorphic vector bundles $\E$ over $\C$ such that ${\rm rank}\,\E=n$, ${\rm deg}\,\E=ng$, and the fiber $\E_\infty$ is spanned by global holomorphic sections. In particular, the locus of $p\in\CC$ such that the fiber $\E_p$ is spanned by global holomorphic sections is a dense set in $\CC$ containing $\infty$.

One may cover a more general moduli space letting the point $\infty$ vary, but for our purposes it is convenient to fix it once for all.

A {\it (global) framing} of $\E\in\moduli$ is a collection $\framing=(\ss_1,\dots,\ss_n)$ of global holomorphic sections $\ss_i\in H^0(\E)$, such that $\ss_1(\infty),\dots,\ss_n(\infty)$ are linearly independent. Let $\wh\moduli$ the set of isomorphism classes of pairs $\left(\E,\framing\right)$ with $\E\in\moduli$ and $\framing=(\ss_1,\dots,\ss_n)$ a framing of $\E$; we require isomorphisms to commute with the framing sections $\ss_i$.

For any $(\E,\framing)\in\wh\moduli$ let us introduce the {\it Tyurin divisor} $\scr T$ of points of $\CC$ where $\ss_1\wedge \dots\wedge \ss_n=0$. By our assumptions on the degree of $\E$, $\mathscr T=\sum_t n_tt$ is a positive divisor of degree $\sum_t n_t=ng$ on $\CC$; moreover, $\infty\not\in\scr T$.

We also fix a conformal disk  $\DD$ such that $\scr T\subset\DD\subset\CC\setminus\{\infty\}$. Let us denote by $z$ the coordinate on $\DD$, writing $z\in\DD$ both for the point and the coordinate representing it.

Given $(\E,\framing)\in\wh\moduli$ we can trivialize the bundle $\E$ over $\CC\setminus\scr T$, where we identify the sections $\ss_i$ with the canonical basis $\mathbf{e}_i$ of $\C^n$. Throughout this paper, all vectors are meant as column vectors unless otherwise explicitly stated.
Denote ${\bf g}_i(z)$ the vectors representing a trivialization of ${\bf s}_i$ on $\DD$ and set 
\be
\label{transitionG}
G(z)=\left[\begin{array}{c|c|c}{\bf g}_1(z)&\cdots&{\bf g}_n(z)\end{array}\right]^t,\qquad z\in\DD.
\ee
The matrix $G(z)$ is a holomorphic function of $z\in\DD$ and $\det G$ vanishes precisely at $\scr T$. The matrix $G$ serves the role of transition function for $\E$; namely, we can identify a local holomorphic section of $\E$ over the open $U\subseteq\CC$ with a vector meromorphic function ${\bf s}={\bf s}(p)$, $p\in U$ with poles at $\scr T\cap U$ only such that
\be
\label{ofirst}
\ss^t(z)G(z)\mbox{ is analytic at }z\in\scr T\cap U.
\ee
For instance, the framing global sections $\ss_i$ are identified with the constant functions ${\bf e}_i$.

Let $\mathcal{R}$ be the ring of holomorphic functions of the variable $z$ in a domain containing the closure $\overline{\mathbb D}$.

\begin{lemma}
\label{lemma1}
Let $G(z)\in\Mat_n(\mathcal{R})$ be such that $\det G(z)\in\mathcal{R}$ does not vanish identically.
There exist unique matrices $H\in\GL_n(\mathcal{R})$ and $\nf\in\Mat_n(\C[z])$ such that 
\be
G(z)= \nf(z) H(z)
\ee
with $\nf$ (termed {\rm polynomial normal form} of $G$) a polynomial matrix of the form 
\be
\label{normform}
\nf(z) = 
\le[
\begin{array}{ccccc}
p_1(z) & 0  & 0 & \cdots & 0\\
f_{21}(z) &p_2(z)  & 0 & \cdots & 0\\
f_{31}(z) & f_{32}(z) &p_3(z)& \cdots & 0\\
\vdots & \vdots & \vdots &\ddots & \vdots \\
f_{n1}(z) & f_{n2}(z) & f_{n3}(z) &\dots &p_n(z) 
\end{array}
\ri]
\ee
where $p_j(z)$ are monic polynomials of degree $d_j$ with all their zeros in $\DD$, and the entries $f_{jk}(z)$ are polynomials of degree $\leq d_j-1$.
Moreover, two matrices $G_1,G_2\in\Mat_n(\mathcal R)$ have the same polynomial normal form if and only if $G_1=G_2 H$ for some $H\in \GL_n(\mathcal R)$. 
\end{lemma}

\noindent {\bf Proof.}
This is a minor elaboration of a similar statement in \cite[Page 137]{Gantmacher1}.
Consider the first row of $G$ and assume it is not identically zero. We find an element whose divisor of zeroes in $\DD$ is of minimal degree and bring it to the $(1,1)$ entry by a column permutation. Now $G_{1,1}(z) = p(z) a(z)$ with $a(z)\in \GL_1(\mathcal{R})$ and $p(z)$ a polynomial with all its zeros in $\DD$. So we can divide the first column by $a(z)$ and assume $G_{1,1}$ is a polynomial. 
By adding a $\mathcal{R}$-multiple of the first column to the remaining columns, we can assume that the whole first row of $G$ is made of polynomials with  $\deg G_{1,j} \leq \deg G_{1,1}-1$. Suppose $G_{1,j}$ are not all zero ($j=2,\dots,n$). We select the one whose divisor in $\DD$ is of  lowest degree and then bring it to the first entry by a column permutation and repeat. 
At each iteration  the degree of $G_{1,1}$  decreases by at least one and the process stops when $G_{1,j},\ j\geq 2$ are zero, while $G_{1,1}$ is a nonzero polynomial (possibly a constant).

We then repeat the process on the submatrix $(2..n; 2..n)$ and so on. At the end of this first pass the matrix $G$ has been transformed to a lower triangular matrix with polynomials on the diagonal whose roots are all in $\DD$.
\be
G \simeq \le[
\begin{array}{cccc}
p_1(z) & 0 & \dots & 0\\
\star  & p_2(z) & \dots & 0\\
 \vdots & &\ddots \\
 \star & \star & \dots  & p_n(z)
\end{array}
\ri], \qquad \deg p_j = d_j.
\ee
The elements on the lower triangular part are some analytic functions, and $\simeq$ denotes identity up to right multiplication by $\GL_n(\mathcal R)$. Again, by column operations only, we can reduce the $(j,k)$ entry to a polynomial of degree $\leq d_j-1$. This proves existence of $\nf,H$ as in the statement. 

To prove  the uniqueness of $\nf,\ H$ it is enough to show that $H=\1$ is the only solution of the equation $\nf=\wt \nf H$ with $\nf,\wt \nf\in\Mat_n(\C[z])$ of the form \eqref{normalization} and $H\in\GL_n(\mathcal R)$. To see it, assume $\nf=\wt \nf H$ and call $h_{ab}$ the entries of $H$. It is clear that $H$ is also lower-triangular and hence for the diagonal entries we have $p_i=\wt p_i h_{ii}$, which is only possible if $h_{ii}=1$ and $p_i=\wt p_i$. Then  $\wt f_{a+1,a}=f_{a+1,a}+h_{a+1,a}p_{a+1}$, and since $\deg (f_{a+1,a}-\wt f_{a+1,a})<\deg p_{a+1}$ we must have $h_{a+1,a}=0$. We proceed to the next sub-diagonal entries, and we have $\wt f_{a+2,a}=f_{a+2,a}+h_{a+2,a}p_{a+2}$, and since $\deg (f_{a+2,a}-\wt f_{a+2,a})<\deg p_{a+2}$ we must have $h_{a+2,a}=0$. Repeating this argument iteratively we prove that $H=\1$. 

For the last statement, suppose $G_i = \nf_i H_i $ for some matrices $\nf_i, H_i$ ($i=1,2$) of the form in the statement.  
If $\nf_1=\nf_2$ then clearly $G_1=G_2 H$ for $H=H^{-1}_2 H_1\in\GL_n(\mathcal R)$. Conversely, if $G_1 = G_2 H$ and $G_1 = \nf H_1$, then $G_2 = \nf H_1 H^{-1}$, and $G_1,G_2$ have the same normal form $\nf$.
\QED

Now, let $G\in\Mat_n(\mathcal R)$ be given in \eqref{transitionG}: we call the coefficients of the polynomials of  $\nf$ in \eqref{normform} the {\it moduli} of the framed bundle $(\E,\framing)\in\wh\moduli$.

To see that the moduli provide local coordinates on $\wh\moduli$, consider the following construction of representative of the equivalence class of framed vector bundles.

\begin{remark}[Construction of a framed bundle from its moduli.]
\rm
\label{remconstruction}
Suppose we have fixed a coordinate disk $\DD\subset\CC$ as above, and within it we fix a polynomial normal form $\nf\in\Mat_n(\C[z])$ as in \eqref{normform} of Lemma \ref{lemma1}, such that the divisor $\scr T:=(\det \nf)$ has degree $ng$; restricting $\nf:\DD\setminus\scr T \to \GL_n(\C)$ we obtain a transition function defining a vector bundle $\E$, trivialized over $\DD$ and $\CC\setminus\scr T$. Namely, as above, local holomorphic sections of the bundle $\E$ just constructed are vector meromorphic functions ${\bf s}={\bf s}(p)$ of $p\in U$ with poles at $\scr T\cap U$ only such that 
\be
\ss^t(z)\nf(z)\mbox{ is analytic at }z\in\mathscr T\cap U.
\ee
We also canonically associate the framing $\framing=(\ss_1,\dots,\ss_n)$ with $\ss_i$ equal to the constant vector ${\bf e}_i$.
It follows that if $\nf\in\Mat_n(\C[z])$ is the polynomial normal of the framed vector bundle $(\E,\framing)$, then the new vector bundle produced by this construction is isomorphic to $\E$. Indeed from the construction of Lemma \ref{lemma1} we have $G=\nf H$ with $H$ analytic and invertible in $\DD$, so that the transition functions $G,\nf$ define isomorphic bundles.
\hfill $\triangle$\end{remark}

Denoting $d_j=\deg p_j$ as in Lemma \ref{lemma1}, we have $\sum_{j=1}^n d_j=ng$. By \eqref{normform} there are $\sum_{j=1}^n jd_j\leq n^2g$ moduli. Therefore, the set $\mathscr M_{\DD}$ of moduli (for a fixed coordinate disk $\DD$) is stratified as
\be
\mathscr M_{\DD}=\bigsqcup\mathscr M_{\DD}^{(d_1,\dots,d_n)},
\ee
where the strata are labeled by $(d_1,\dots,d_n)$ with $d_j\geq 0,\sum_{j=1}^nd_j=ng$; the moduli provide a complex manifold structure of dimension $\sum_{j=1}^njd_j$ to $\mathscr M^{(d_1,\dots,d_n)}_\DD$. The open stratum $\mathscr M^{(0,\dots,0,ng)}_\DD$, of dimension $n^2g$, is described in the following example.

\begin{example}\label{example:generic}
According to Lemma \ref{lemma1}, in the {\rm generic} case $(d_1=0,\dots,d_n=ng)$ the matrix $G$ has the form
\be
G(z)= \le[
\begin{array}{ccccc}
1 &0& \cdots & 0 &0\\
0 & 1 & \cdots & 0 &0\\
\vdots & \vdots & \ddots & \vdots & \vdots \\
0 & 0 & \cdots & 1 & 0\\
h_{1}(z) &h_2(z) &\cdots &h_{n-1}(z) & p(z)
\end{array}
\ri] H(z)
\ee
where $H(z)\in\GL_n(\mathcal{R})$, $p(z)$ is a monic polynomial of degree $ng$ with zeroes in $\DD$ and $h_j(z)$ are polynomials in $\C[z] /\langle p(z)\rangle$, i.ee., arbitrary polynomials of degree $\leq ng-1$. The number of parameters of the polynomials $(h_j)_{j=1}^{n-1}$ is $(n-1)ng$, thus providing, along with the coefficients of $p$, a total number of $n^2g$ parameters.
In geometric terms this shows that $\E$ is generically an extension of a rank $n-1$ trivial subbundle by the line bundle $\det\E$;
\be
\nn
0\to\O^{n-1}\to\E\to\det\E\to 0.
\ee
\end{example}

For different coordinate disks $\DD\subset\CC$ this construction provides an atlas of complex charts of dimension $n^2g$ on the open part in $\wh\moduli$ where the normal form is generic in the sense of this example. 

\br
Lemma \ref{lemma1} also shows that the transition function can always be put in lower-triangular form. Hence there is a complete flag $\mathcal{E}_0=\left\lbrace 0\right\rbrace\subset\mathcal E_1\subset\cdots\subset \E_n=\E$ of sub-bundles of $\E$, with $\mathcal E_r$ of rank $r$ (see e.g. \cite{Gunning}). Clearly $\deg\E_r=\sum_{i=1}^rd_i$; in particular (semi-)stability  implies $\sum_{j=1}^r d_j < r g$ ($\sum_{j=1}^r d_j \leq r g$) for all $r\in\{1,\dots, n\}$.
\er

\subsection{Tyurin space}

Let us denote, as above, $\mathcal R$ the ring of functions of $z$ holomorphic in a domain containing $\overline{\DD}$; let also $\vec{\mathcal R}=\C^n\otimes\mathcal R$, thought of as row-vectors.
Moreover, for a (possibly vector-valued) meromorphic function $f=f(z)$ of $z\in\DD$, with finitely many poles, we define
\be
\label{Cauchyoperator}
{\sf C}_-[f](z):=
\sum_{p\in(f)_-}\res{z'=z(p)}\frac{f(z')\d z'}{z-z'},
\ee
where the sum extends over the poles of $f$ in $\DD$, the Cauchy operator extracting the polar part of $f$ in $\DD$.

\begin{definition}
The {\rm Tyurin space} associated with the framed bundle $(\E,\framing)\in\wh\moduli$ is the complex vector space
\be
\label{tyurinspace}
T:={\sf C}_-\left[\vec{\mathcal R}\nf^{-1}\right],
\ee
where $\nf$ is the polynomial normal form introduced above.
\end{definition}

Let us introduce the vectors ${\bf v}_1^t(z),\dots,{\bf v}_{ng}^t(z)\in T$ defined by
\be
\label{tyurinbasis}
{\bf v}_{d_1+\cdots+d_{j-1}+k}^t(z)={\sf C}_-[z^k{\bf e}_j^t P^{-1}(z)], \qquad k=0,\dots,d_j-1,\ j=1,\dots,n,
\ee
where $d_j$ is the degree of the diagonal entries $p_j(z)$ in the polynomial normal form $\nf$, as in \eqref{normform} (we recall that $\sum_{j=1}^n d_j=ng$).
The proof of the following lemma is elementary and hence omitted.

\begin{lemma}
The vectors in \eqref{tyurinbasis} form a basis of $T$; in particular $\dim T=ng$.
\end{lemma}

\begin{remark}[Relationship with Tyurin vectors]
\label{rem:simpleT}
\rm
If the Tyurin divisor $\scr T$ associated with $(\E,\framing)\in\wh\moduli$ is denoted $\scr T=\sum_j n_jt_j$, then $T$ is a direct sum of spaces of Laurent series in $z-z(t_j)$ of dimensions $n_j$; moreover ${\sf C}_-[z T]\subseteq T$, and in particular each of the ``localized" subspaces is invariant under multiplication by $z-z(t_j)$ and projection to the polar tail. If in addition, the Tyurin divisor is simple (i.e. $n_j=1$) we have
\be
T=\bigoplus_{j=1}^{ng}\C{\bf v}^t_j,\qquad \mathbf{v}_j(z):=\frac{{\bf h}_j}{z-z(t_j)}\label{TyurinVectors}
\ee
for some nonzero vectors ${\bf h}_1,\dots,{\bf h}_{ng}\in\C^n$ ({\it Tyurin vectors} \cite{EnRub,krichever,Tyurin1, Tyurin2}).
\hfill $\triangle$\end{remark}

\begin{remark}[Cohomological interpretation of the Tyurin space]\label{remcoh}
\rm
Given $(\E,\framing)\in\wh\moduli$, we can define a short exact sequence of sheaves on $\CC$
\be
\label{sesTyurin}
0\to\O^n\to\E\to\mathcal{T}\to 0.
\ee
Here $\O^n$ is the sheaf of $\C^n$-valued holomorphic functions on $\CC$, and $\E$ is identified with the sheaf of its holomorphic sections; the first map is defined, for any open $U\subseteq\CC$, by
\be
(f_1,\dots,f_n)\in\O^n(U)\mapsto f_1\ss_1+\cdots+f_n\ss_n\in\E(U).
\ee
The quotient sheaf $\mathcal T$ is supported at the Tyurin divisor $\mathscr T$. 
The Tyurin space $T$ introduced above in \eqref{tyurinspace} is isomorphic to the space $H^0(\cal T)$ of global sections of $\mathcal T$; indeed, $H^0(\mathcal T)$ can be regarded as the space of polar tails of rational functions of $z\in\DD$, such that right-multiplied by $\nf(z)$ they become analytic on $\DD$.
\hfill $\triangle$\end{remark}

\section{Non-abelian theta divisor}
\label{BrillNoether}

Let $\thetadiv\subset\moduli$ the set of isomorphism classes of vector bundles $\E$ with $h^0(\E)>n$. Let us note that, in view of the Riemann--Roch theorem and of the Serre duality, this condition is equivalent to $h^1(\E)=h^0(\E^*\otimes \mathcal K)>0$.
Similarly, let $\wh{\thetadiv}\subset\wh\moduli$ consist of framed bundles $(\E,\framing)$ with $h^0(\E)>n$.
Let us also denote $\moduli_0=\moduli\setminus\thetadiv$ and $\wh\moduli_0=\wh\moduli\setminus\wh\thetadiv$.
We shall refer to $\thetadiv$ and $\wh{\thetadiv}$ as {\it non-abelian theta divisors}; the goal of this section is to provide explicit holomorphic equations for $\wh{\thetadiv}$.

Let $T$ be the Tyurin space \eqref{tyurinspace} associated with $(\E,\framing)\in\wh\moduli$.
Denote $\mathcal K$ the canonical bundle on $\CC$ and introduce a bilinear pairing between $T$ and $H^0(\mathcal K^n)$ as follows;
\be
\label{bilinear}
\mathbf{v}^t\otimes\bs\nu\mapsto\left\langle\mathbf{v}^t,\bs\nu\right\rangle=
\sum_{t\in\scr T}\res{t}\mathbf v^t\bs\nu.
\ee
Hereafter we agree that elements of $H^0(\mathcal K^n)$ are (column) vector holomorphic differentials on $\CC$.

Let us fix a basis $\omega_1,\dots,\omega_g$ of holomorphic differentials on $\CC$ and let ${\bf e}_j$ be the standard basis of $\C^n$. With respect to the basis ${\bf v}_j^t(z)$ \eqref{tyurinbasis} of $T$, and the basis $\omega_i\otimes{\bf e}_j$ of $H^0(\mathcal K^n)$, the pairing \eqref{bilinear} is represented by the {\it Brill--Noether--Tyurin matrix}
\be
\label{matrix}
\mathbb T=\sum_{t\in\scr T}\res{z=z(t)}\left[\begin{array}{c}
{\bf v}_1^t(z)\\
\vdots\\
{\bf v}_{ng}^t(z)
\end{array}\right]\otimes\left[\omega_1(z),\dots,\omega_g(z)\right]\in\Mat_{ng\times ng}(\C),
\ee
where $\otimes$ denotes the Kronecker product.
Note that $\mathbb T$ depends holomorphically on the moduli of $(\E,\framing)$.

We explicitly identify the left and right null-spaces of this pairing, i.e., respectively, the kernels of the maps $T\to H^0(\mathcal K^n)^*$ and $H^0(\mathcal K^n)\to T^*$ naturally induced by the bilinear pairing $\langle,\rangle$.
To this end, given $(\E,\framing)\in\wh\moduli$ we consider the bundle $\E^*\otimes\mathcal K$ trivialized similarly as explained above for $\E$; namely we identify its local holomorphic sections over the open $U\subseteq\CC$ with vector holomorphic differentials $\bs\eta=\bs\eta(p)$ of $p\in U$ satisfying
\be
\label{rhpdualvector}
\nf^{-1}(z)\bs\eta(z)\mbox{ analytic at }z\in\scr T\cap U,
\ee
with $\nf(z)$ the polynomial normal form.

\begin{proposition}
\label{propdual}
\noindent{\bf [1]} The left null-space of \eqref{bilinear} is identified with the space $H^0(\E\otimes\O(-\infty))$ of global holomorphic sections of $\E$ vanishing at $\infty$ by the map
\be
\label{map1}
{\bf r}\in H^0(\E(-\infty))\mapsto {\sf C}_-[{\bf r}(z)]\in T.
\ee
In particular 
\be
\label{corank1}
h^0(\E)=n+{\rm corank}\,\mathbb T.
\ee
\noindent{\bf [2]} The right null-space of the pairing \eqref{bilinear} is identified with $H^0(\E^*\otimes \mathcal K)$. 
In particular
\be
\label{corank2}
h^0(\E^*\otimes \mathcal K) = h^1(\E) = {\rm corank}\, \mathbb T.
\ee 
\end{proposition}

\noindent{\bf Proof. }
{\bf [1]} The section ${\bf r}\in H^0(\E(-\infty))$ by the above discussion is a vector meromorphic function ${\bf r}={\bf r}(p)$ of $p\in\CC$ with ${\bf r}(\infty)=0$ and poles at $\scr T$ only, such that ${\bf r}^t\nf$ is analytic at $\scr T$. For all vectors of holomorphic differentials $\bs\nu\in H^0(\mathcal K^n)$ we have
\be
\langle {\sf C}_-[{\bf r}(z)],\bs\nu\rangle=\sum_{t\in\scr T}\res{t}{\bf r}^t\bs\nu=0,
\ee
as ${\bf r}^t\bs\nu$ is a globally defined meromorphic differential on $\CC$ with poles at $\scr T$ only.

Conversely, given ${\bf v}^t(z)\in T$ and fixing a symplectic basis $\{a_1,b_1,\dots,a_g,b_g\}$ of $H_1(\CC,\mathbb Z)$, there exists a unique row vector meromorphic differential $\bs\Omega^t$ on $\CC$ satisfying
\be
(\bs\Omega)\geq -\scr T,\qquad \bs\Omega(z)-\d{\bf v}^t(z)\mbox{ analytic in }\DD,\qquad\oint_{a_j}\bs\Omega=0,\ j=1,\dots,g.
\ee
Define a row vector function ${\bf r}$ by ${\bf r}(p)=\int_\infty^p\bs\Omega$; ${\bf r}(p)$ is a single-valued vector meromorphic function of $p\in\CC$ if and only if $\bs\Omega$ has vanishing $b$-periods too.
It follows from the Riemann bilinear relations for second kind abelian differentials that this is the case if and only if ${\bf v}^t(z)\in T$ is in the left null-space of $\langle,\rangle$.
In such case, $\bf r\in H^0(\E(-\infty))$ and we have explicitly inverted \eqref{map1}.

\noindent {\bf [2]} It follows from the fact that $\nf^{-1}\bs\eta$ is regular at $\scr T$ if and only if ${\bs\eta}$ is in the right null-space of $\langle,\rangle$.
\QED

Incidentally, the Riemann--Roch formula \eqref{riemannroch} also follows by comparing \eqref{corank1} and \eqref{corank2}.

\begin{corollary}
The locus $\wh{\thetadiv}$ is a divisor in $\wh\moduli$, defined by local holomorphic equations $\det\mathbb{T}=0$.
\end{corollary}

The fact that $\det\mathbb T$ is not identically zero on $\wh\moduli$ follows by considering the bundle $\E=\O(\mathscr D)^{\oplus n}$ for a non-special degree $g$ divisor $\scr D$ on $\CC$; taking $\scr D$ special shows instead that $\wh\thetadiv$ is non-empty.
 
\begin{example}
\label{example2}
In the case of line bundles, $n=1$, and distinct Tyurin points $t_1,\dots,t_g$, $\mathbb{T}$ reduces to the Brill--Noether matrix
$\mathbb T = \left[\omega_i(t_j)\right]_{i,j=1}^g$ and the condition $\det\mathbb{T}=0$ defines the locus of degree $g$ special divisors.
\end{example} 

\begin{remark}[Cohomological interpretation of the pairing]
\rm
The short exact sequence \eqref{sesTyurin} induces the long exact sequence in cohomology
\be
0\to H^0(\O^n)=\C^n\to H^0(\mathcal{E})\to H^0(\mathcal{T})\to H^1(\O^n)\to H^1(\mathcal{E})\to 0
\ee
where we note that $H^1(\mathcal T)=0$ because $\mathcal T$ is a flasque sheaf (i.e., restriction morphisms are surjective) hence its $i$th cohomology groups vanish for $i\geq 1$~\cite[Proposition~2.5, Chapter~III]{Hartshorne}).
By Serre duality and the identification of $H^0(\mathcal{T})$ with the Tyurin space $T$ (Remark~\ref{remcoh}), we rearrange this exact sequence as
\be
\label{lescoh}
0\to H^0(\E(-\infty))\to T\xrightarrow{\phi} H^0(\mathcal K^n)^*\to H^0(\mathcal{E}^*\otimes\mathcal K)^*\to 0,
\ee
where we identify $H^0(\E)/\C^n\simeq H^0(\E(-\infty))$.
It follows that the map $\phi$ induces a pairing 
\be
\langle\ ,\ \rangle:T \otimes H^0(\mathcal K^n)\to\mathbb{C}.
\ee
By exactness of \eqref{lescoh} the left null-space is isomorphic to $H^0(\mathcal{E}(-\infty))$ and right null-space is isomorphic to $H^0(\mathcal{E}^*\otimes\mathcal K)$, as shown explicitly in Proposition~\ref{propdual}.
\hfill $\triangle$\end{remark}

\begin{remark}[Non-abelian theta divisor and semi-stability]
\label{remstable}
\rm
If $\E\in\moduli_0$ then $\E$ is necessarily semi-stable. For, suppose $\E$ is not semi-stable and let $\mathcal F$ be a sub-bundle of rank $r\not=0,n$ such that $\deg\mathcal F=rg+k$ for some $k\geq 1$. On one side, the Riemann--Roch Theorem applied to $\mathcal F$ implies $h^0(\mathcal F)=h^1(\mathcal F)+\deg\mathcal F-r(g-1)\geq r+k$; on the other side, since $\E$ is generically spanned, we must have at least another $n-r$ sections of $\E$, and so $h^0(\E)\geq n+k$.\footnote{We thank Prof. Indranil Biswas for pointing this out.}
\hfill $\triangle$\end{remark}

\section{Non-abelian Cauchy kernel}
\label{secCauchy}

We now introduce the main tool in the subsequent analysis, namely the matrix Cauchy kernel. The construction is explicit in terms of the moduli of the bundle. 
\begin{proposition}
\label{propdualmero}
Let $(\E,\framing)\in\wh\moduli_0$. Fix a point ${p} \in\CC\setminus(\scr T+\infty)$. There exist exactly $n$ linearly independent global holomorphic sections of the bundle $\E^*\otimes\mathcal K(p+\infty)$.
\end{proposition}

The statement of the proposition is equivalent to $h^0(\E^*\otimes\mathcal K(p+\infty))=n$, i.e., by Riemann--Roch, to $h^0(\E(-p-\infty))=0$; however, the assumptions imply that no nontrivial global holomorphic section of $\E$ can vanish at $\infty$, proving the proposition.
For later convenience we give a direct proof which does not make use of Riemann--Roch's theorem and which instead uses the Tyurin data.
\\

\noindent{\bf Proof.}
A global holomorphic section of $\E^*\otimes\mathcal K(p+\infty)$ is a vector meromorphic differential $\bs\nu$ such that $\bs\nu$ has simple poles at $p,\infty$ only and $\nf^{-1}\bs\nu$ is analytic at $\scr T$.
Let us fix a third kind abelian differential $\omega_{{p},\infty}$ on $\CC$ with simple poles at ${p},\infty$ only; $\boldsymbol{\nu}$ must be in the form $\boldsymbol{\nu}=\omega_{{p},\infty}\mathbf{e}+\wh{\boldsymbol{\nu}}$, where $\wh{\boldsymbol{\nu}}$ is a vector of holomorphic differentials on $\CC$ and $\bf e\in\C^n$ is any constant vector.
The condition that $\nf^{-1}\boldsymbol{\nu}$ does not have poles at $\scr{T}$ translates into the condition
\be
\sum_{t\in\scr T}\res{z=z(t)}{\bf v}^t(z)\left[\omega_{{p},\infty}(z)\mathbf{e}+\wh{\boldsymbol{\nu}}(z)\right]=0,\qquad\forall{\bf v}^t(z)\in T.
\label{conditiontwopoles}
\ee
The condition uniquely determines $\wh{\bs\nu}$ due to the assumption $\E\not\in\thetadiv$. Indeed, recalling the basis\eqref{tyurinbasis} of $T$, fixing a basis $\omega_1,\dots,\omega_g$ of $H^0(\mathcal K)$, and writing $\wh{\boldsymbol{\nu}}=\sum_{k=1}^g\mathbf{b}_k\omega_k$ for the unknown constant vectors $\mathbf{b}_1,\dots,\mathbf{b}_g\in\mathbb{C}^n$, the condition \eqref{conditiontwopoles} is equivalent to
\be
\label{linearinho}
\mathbb{T}
\renewcommand*{\arraystretch}{1.2}
\left[\begin{array}{c}
\mathbf{b}_1 \\ \hline \vdots \\ \hline \mathbf{b}_g
\end{array}\right] + 
\sum_{t\in\scr T}\res{z=z(t)}\left[
\begin{array}{c}
\omega_{{p},\infty}(z){\bf v}_1^t(z)
\\ \hline
\vdots
\\ \hline
\omega_{{p},\infty}(z){\bf v}_{ng}(z)
\end{array}
\right]\mathbf{e}=0,
\ee
where $\mathbb T$ is given in \eqref{matrix}, for the same choice of $\omega_1,\dots,\omega_g$. This has a unique solution $\mathbf{b}_1,\dots,\mathbf{b}_g$, as $\det\mathbb{T}\not=0$ since $\E\not\in\thetadiv$. The conclusion follows as $\mathbf{e}$ is arbitrary in $\mathbb{C}^n$.
\QED

\begin{shaded}
\begin{theoremdefinition}[Non-abelian Cauchy kernel]
\label{defCauchy}
Fix a framed vector bundle $(\E,\framing)\in \wh\moduli_0$. There exists a unique $n\times n$ {\rm (matrix) Cauchy kernel} $\Cauchy_\infty({q},{p})$ such that 
\begin{itemize}
\item[-] it is a meromorphic differential in the variable ${q}$ with its divisor satisfying  $(\Cauchy_\infty({q},{p}))_{q}\geq -{p}-\infty$;
\item[-] it is a meromorphic function in the variable ${p}$ with its divisor satisfying $(\Cauchy_\infty({q},{p}))_{p}\geq \infty-{q}-\scr{T}$,
\end{itemize}
and such that  the following properties hold:
\begin{enumerate}
\item the expression $\nf^{-1}(q)\Cauchy_\infty({q},{p})$ is regular for ${q}\in\DD$, and
\item it satisfies the normalization condition
\be
\label{normalization}
\res{{q}={p}}\Cauchy_\infty({q},{p})=\1_{n\times n}=-\res{{q}=\infty}\Cauchy_\infty({q},{p}).
\ee
\end{enumerate}
\end{theoremdefinition}
\end{shaded}
 
\noindent{\bf Proof.}
According to Proposition \ref{propdualmero}, let $\boldsymbol{\eta}_{\ell}(q;{p},\infty)$ be the unique sections of $\E^*\otimes\mathcal K$ (in $q$) with simple poles at $q={p},\infty$ only, normalized by 
\be\label{res}
\res{{q}={p}}\bs\eta_{\ell}({q};p,\infty)=\mathbf{e}_\ell=-\res{{q}=\infty}\boldsymbol{\eta}_\ell({q};p,\infty),\qquad\ell=1,\dots,n.
\ee
Then the matrix
\be
\label{construction}
\Cauchy_\infty({q},{p}):=\left[
\begin{array}{c|c|c}
\boldsymbol{\eta}_1(q;{p},\infty) & \cdots & \boldsymbol{\eta}_n(q;{p},\infty)
\end{array}
\right]
\ee
has the desired properties. The uniqueness follows from Proposition~\ref{propdualmero}.
\QED

We can write the Cauchy kernel more explicitly in terms of the basis ${\bf v}_j^t(z)$ given in \eqref{tyurinbasis}, of a fixed basis $\vec{\omega}=(\omega_1,\dots,\omega_g)$ (understood as a row-vector) of holomorphic differentials on $\CC$, and of the matrix $\mathbb{T}$ given in \eqref{matrix}.
Indeed from \eqref{linearinho} and \eqref{construction} we infer that the Cauchy kernel can be written as
\be
\label{expressioncauchyblock}
\Cauchy_\infty({q},{p})=\omega_{{p},\infty}({q})\1_{n\times n} -  (\vec\omega({q})\otimes\1_{n\times n})\cdot\mathbb{T}^{-1}\cdot
\sum_{t\in\scr T}\res{z=z(t)}
\left[
\renewcommand*{\arraystretch}{1.1}\begin{array}{c}
\mathbf{v}_1^t(z)
\\ \hline
\vdots
\\ \hline\mathbf{v}_{ng}^t(z)
\end{array}
\right]\omega_{{p},\infty}(z)
\ee
where $\otimes$ denotes the Kronecker product, and $\omega_{{p}_+,{p}_-}({q})$ denotes the third kind differential with simple poles at ${p}_\pm$ and residue $\pm 1$, respectively, normalized to have vanishing $a$-cycles. In more transparent terms, the entries of $\Cauchy_\infty({q},{p})$ can be expressed as
\begin{align}
\label{expressioncauchydet}
\left(\Cauchy_\infty({q},{p})\right)_{i,j}&=\frac{\det
\renewcommand*{\arraystretch}{1.4}
\left[\begin{array}{c|c}
\mathbb{T} & \sum_{t\in\scr T}\res{z=z(t)}
\left[\begin{array}{c}
\mathbf{v}_1^t(z)\mathbf{e}_i\omega_{{p},\infty}(z) \\ 
\hline \vdots \\ \hline
\mathbf{v}_{ng}^t(z)\mathbf{e}_i\omega_{{p},\infty}(z)
\end{array}\right]
\\ \hline
\vec\omega({q})\otimes \mathbf e^t_j & \delta_{ij}\omega_{{p},\infty}({q})
\end{array}\right]
}{\det \mathbb{T}}.
\end{align}
From the expression \eqref{expressioncauchydet} it follows that $\Cauchy_\infty({p},{q})$ is a {\it meromorphic} function of the moduli with pole at $\wh{\thetadiv}$ only.
 
It follows from the uniqueness of the Cauchy kernel, that \eqref{expressioncauchydet} is single-valued in ${p}$. This can also be checked directly from the monodromy properties of normalized third kind differentials
\be
\omega_{{p}+{\a_\ell},\infty}({q})=\omega_{{p},\infty}({q}), \qquad
\omega_{{p}+{\b_\ell},\infty}({q})=\omega_{{p},\infty}({q})+\omega_\ell({q})
\ee
with respect to analytic continuations along a symplectic basis of generators $\a_1,\b_1,\dots,\a_g,\b_g$ for $\pi_1(\CC,\infty)$ for which $\oint_{\a_\ell}\omega_k=\delta_{\ell k}$; from this we see that the determinant of the matrix in the numerator of \eqref{expressioncauchydet} is single-valued.

Furthermore, \eqref{expressioncauchydet} shows that the Cauchy kernel has poles as a function of ${p}$ as ${p}\to \scr{T}$. To see this, for $p\in\DD$ we have
\be
\label{divergence}
\sum_{t\in\scr T}\res{z=z(t)}\omega_{{p},\infty}(z)\mathbf{v}_j^t(z)=-\res{z=z(p)}\omega_{{p},\infty}(z)\mathbf{v}_j^t(z)=-\mathbf{v}_j^t(z(p))
,\quad j=1,\dots,ng,
\ee
which diverges as $p\to\scr T$. However, \eqref{divergence} shows that the singular part of $\Cauchy_\infty(q,p)$ as ${p}\to\scr T$ belongs to the Tyurin space $T$ (row-wise); then $\Cauchy_\infty(q,p)\nf(p)$ is regular as $p\to\scr T$. 

\begin{corollary}
\label{corollCauchy}
The expression $\nf^{-1}(q)\Cauchy_\infty(q,p)\nf(p)$ is regular for  $q,p\in \DD$ away from the diagonal locus $q=p$.
\end{corollary}

\begin{example}In the case of line bundles, $n=1$, $\Cauchy_\infty({q},{p})$ reduces to the known (scalar) Cauchy kernel \cite{Fay73}; e.g. when the Tyurin points $t_1,\dots,t_g$ are simple, we have
\be
\Cauchy_\infty({q},{p})=\frac{\det\left[\begin{array}{cccc}
\omega_1(t_1) & \cdots & \omega_g(t_1) & \omega_{{p},\infty}(t_1) \\
\vdots & \ddots & \vdots & \vdots \\
\omega_1(t_g) & \cdots & \omega_g(t_g) & \omega_{{p},\infty}(t_g) \\
\omega_1({q}) & \cdots & \omega_g({q}) & \omega_{{p},\infty}({q})
\end{array}\right]}{\det\left[\begin{array}{ccc}
\omega_1(t_1) & \cdots & \omega_g(t_1) \\
\vdots & \ddots & \vdots \\
\omega_1(t_g) & \cdots & \omega_g(t_g)
\end{array}\right]}.
\ee
\end{example}

\begin{example}[Cauchy kernel for simple Tyurin points]
\label{exTyusymp}
Consider the case of simple Tyurin point, i.e., $\scr T=t_1+\cdots+t_{ng}$ with $t_i\not=t_j$ for $i\not=j$, and let $\bf h_i$ be the Tyurin vectors, cf. Remark \ref{rem:simpleT}.
Define $\mathbf H = [{\bf h}_1,\dots ,{\bf h}_{ng}]\in \Mat_{n, ng}(\C)$ (the normalization of the Tyurin vectors are chosen arbitrarily, and the final formul\ae\ are independent of this choice), and, for any meromorphic differential $\omega$, define
\be
\omega(\scr T):= {\rm diag} (\omega(t_1), \dots, \omega(t_{ng})).
\ee
Then the Brill--Noether--Tyurin matrix is defined by
\be
\label{BNTsimple}
\mathbb T = \le[\omega_1(\scr T) {\bf H}^t|\cdots |\omega_g(\scr T) {\bf H}^t\ri]\in \Mat_{ng\times ng}
\ee
where $\omega_j$ is any basis of holomorphic differentials (for example a normalized basis).
The non-abelian Cauchy kernel is given explicitly in terms of the normalized abelian differentials of the first and third kind~\cite[Section~9.10]{Forster} as the following expression
\bea
\label{Cauchyentries}
\renewcommand*{\arraystretch}{1.2}
\big(\Cauchy_\infty(q,p)\big)_{ij}=\frac{
\det \le[\begin{array}{c|c|c|c}
\begin{array}{c}
\\
\omega_1(\scr T) \mathbf H^t
\\
\phantom{j}\end{array}  & \cdots &\begin{array}{c}
\\
\omega_g(\scr T) \mathbf H^t
\\
\phantom{j}
\end{array}  & 
\omega_{p,\infty} (\scr T) \mathbf H^t \mathbf e_j\\
\hline
 \omega_1(q){\mathbf e}_i^t & \cdots & \omega_g(q) \mathbf e_i^t & \omega_{p,\infty}(q)\delta_{ij}
\end{array}\ri]
}
{\det \mathbb T},\ \ i,j=1,\dots, n.
\eea
Here the non-abelian theta divisor is easily described as the locus $\det\mathbb T=0$  and the corank of this matrix is precisely $h^1(\E)$ (see Proposition~\ref{propdual}).
\end{example}

\section{Framed Higgs fields and \texorpdfstring{$T^*\wh\moduli_0$}{T*V0}}
\label{secHiggs}

Given $\E\in\moduli$, a {\it framed Higgs field} $\Phi$ is an element $\Phi\in H^0(\mathrm{End}(\E)\otimes\mathcal K(\infty))$, namely a meromorphic section of $\mathrm{End}(\E)\otimes\mathcal K$ with only a simple pole at $\infty$.
Given a framing $\framing$ of $\E$, inducing the polynomial normal form $\nf$, a framed Higgs field is trivialized as a matrix meromorphic differential $\Phi$ with poles at $\mathscr T+\infty$ only, such that
\be
\label{boundaryhiggs}
\nf^{-1}(z)\Phi(z)\nf(z)\mbox{ is analytic on }z\in\DD.
\ee

The next lemma, whose proof we omit,  provides an equivalent characterization of framed Higgs fields, and is a generalization of a result of \cite{krichever}.

\begin{lemma}
\label{lemmaframedhiggs}
A meromorphic matrix differential $\Phi$ on $\CC$ with $(\Phi)\geq -\mathscr{T}-\infty$ defines a framed Higgs field if and only if ${\sf C}_-[\Phi]\in T$ (row-wise) and ${\sf C}_-[T\Phi]\subseteq T$.
Here $T$ is the Tyurin space \eqref{tyurinspace}, ${\sf C}_-$ is the Cauchy operator \eqref{Cauchyoperator}, and $T\Phi$ denotes the set of elements ${\bf v}^t\Phi$ for all ${\bf v}^t\in T$.
\end{lemma}

The Cauchy kernel can be used to construct Higgs fields explicitly in terms of local data near the Tyurin points as we now explain.

\begin{proposition}
\label{propHiggs}
Given an arbitrary collection $(\phi_t)_{t\in \scr T}$ of germs of matrix valued holomorphic differentials, then the expression
\be
\label{omegaphi}
\Phi({q}):= \sum_{t\in\scr T}\res{{p}=t}\Cauchy_\infty({q},{p}) \nf(p) \phi_{t}(p)\nf^{-1}(p)
\ee
defines a Higgs field. Viceversa, given any Higgs field $\Phi$, then the expression at the right hand side of \eqref{omegaphi} with $\phi_t(p)$ equal to the germ of $\nf^{-1}(p)\Phi(p)\nf(p)$ at $p=t$ recovers $\Phi$.
\end{proposition}
The proof is a simple verification based on the properties of the Cauchy kernel in Thm-Def. \ref{defCauchy} and it is omitted. 

It is well known that the Higgs fields are generically identifiable with the cotangent vectors to the moduli space of vector bundles; we need to make this identification explicit in our setting.
To this end, {following \cite{krichever}} introduce the following pairing between a framed Higgs field $\Phi\in H^0(\mathrm{End}(\E)\otimes\mathcal K(\infty))$ and a deformation $\pa$ of the transition function $G$:
\be
\label{pairing1}
\left(\Phi,\pa \right)=\sum_{t \in \scr T} \res{t} \tr \left(\Phi \pa G G^{-1} \right).
\ee
Note that this expression is invariant under $G\mapsto G H$ with $H$ analytic and invertible in $\DD$. Therefore we may replace $G$ in \eqref{pairing1} by the normal form $\nf$. 

\noindent { The next Proposition \ref{propnondegenerate} can be found in \cite{krichever} but we provide a short proof for the reader's convenience.}
\begin{proposition}
\label{propnondegenerate}
Let $\E\in \wh \moduli_0$: then the pairing \eqref{pairing1} is nondegenerate.
\end{proposition}

\noindent {\bf Proof.}
We need to show that \eqref{pairing1} vanishes  for all framed Higgs fields $\Phi \in H^0(\mathrm{End}(\E)\otimes\mathcal K(\infty))$ if and only if the deformation is trivial. 
The only nontrivial part is to show that if $(\Phi,\pa)=0$ for all framed Higgs fields then the deformation is trivial.
Let us compute explicitly the pairing when the Higgs field is expressed via \eqref{omegaphi}.

Let $(\phi_t)_{t\in \scr T}$ be a collection of germs of matrix valued holomorphic differentials at $t\in\scr T$, and let $\Phi$ be the Higgs field defined in \eqref{omegaphi}. Then
\be
(\Phi,\partial)=\sum_{t\in\scr T}\res{t}\tr\left(\phi_{t} \nf ^{-1}\pa \nf \right).
\ee
This is seen by inserting the expression \eqref{omegaphi} in \eqref{pairing1} (with $G$ replaced by $P$) and then evaluating the residues using the properties of the Cauchy kernel.
Since now the  collection $(\phi_t)$ of holomorphic matrix germs  at $t\in\scr T$ is completely arbitrary, we must have $\pa \nf =0$.
\QED

\paragraph{From framed to unframed vector bundles.}
There is an action of $\P\GL_n(\C)$ on $\wh\moduli$ defined by a change of global framing; it is induced by the $\GL_n(\C)$-action
of multiplication of the transition function $G$ on the left by a given constant matrix $C\in \GL_n(\C)$.

The orbits of this $\P\GL_n(\C)$-action on $\wh\moduli_0$ are in one-to-one correspondence with points of $\moduli_0$; 
\be\label{quotPGL}
\moduli_0=\wh\moduli_0/\P\GL_n(\C).
\ee

Therefore the cotangent space $T^*_\E\moduli_0$ is identified with the subspace of framed Higgs fields $\Phi$ annihilating of the infinitesimal version $\pa G=AG$, $A\in\Mat_n(\C)$, of the deformation; using  \eqref{pairing1}  we get
\be
\label{mmmm}
(\Phi,\pa)=\sum_{t\in \scr T} \res{t} \tr \left(\Phi\pa G G^{-1}\right)=\sum_{t\in \scr T} \res{t} \tr \left(\Phi A\right)=-\res{\infty}\tr\left(\Phi A\right)=0.
\ee
Since the matrix $A\in\Mat_n(\C)$ in \eqref{mmmm} is arbitrary we immediately get the following well-known result. 

\begin{proposition}
\label{cotang} 
The cotangent space $T^*_\E\moduli_0$ is identified with the space $H^0({\rm End}(\E)\otimes\mathcal K)$ of {\rm holomorphic} Higgs fields.
\end{proposition}

\section{Flat connections and the Cauchy kernel}
\label{secAffine}

The previous sections have described the holomorphic moduli on $\moduli_0$ (and $\wh \moduli_0$) and the identification of the cotangent bundle in terms of the explicitly realized Higgs fields in terms of local holomorphic germs (Proposition~\ref{propHiggs}) together with the tautological one form on the cotangent bundle \eqref{pairing1}.

The main tool is the Cauchy kernel (Thm-Def. \ref{defCauchy}) of the framed bundle which  allows us to define a holomorphic map from the cotangent bundle $T^*\moduli _0$ to the $\GL_n(\C)$ character variety. This hinges on the definition of the ``Fay'' differential as the regular term in the  expansion of $\Cauchy_\infty$ along the diagonal and subsequently by Theorem \ref{nosing}.

\subsection{Affine bundles of connections}
Let us fix a {non-special} positive divisor of  degree $g$,  $\scr D=\sum_p n_pp$ (i.e. $h^0(\mathcal K(-\scr D))=0$) with  $\infty\not\in\scr D$.

For an arbitrary $\E\in \moduli_0$, let $\conn^{\scr D}_\E$ be the affine space of holomorphic connections on $\E^*\otimes \mathcal O (\scr D)$; this is an affine space modelled on $H^0(\E^*\otimes \E\otimes\mathcal K)$ which is isomorphic to  $T^*_\E\moduli_0$  by Proposition~\ref{cotang}.

Let $\conn^{\scr D}$ be the affine bundle over $\moduli_0$, whose fiber at $\E\in\moduli_0$ is $\conn^{\scr D}_\E$; equivalently, this is the set of pairs $(\E,\nabla)$ of $\E\in\moduli_0$ and $\nabla\in\conn^{\scr D}_\E$.

The final goal of this section is to provide a reference connection $\nabla_\scr D$ in the affine bundle  $\conn^{\scr D}\to\moduli_0$ to be  identified with the zero section of $T^*\moduli_0\to\moduli_0$. We also remark that if $\E\in \thetadiv$ the space of connections $\conn^{\scr D}_\E$ is still an affine space modelled over $H^0(\E^*\otimes \E\otimes\mathcal K)$, but it has larger dimension and hence $\conn^{\scr D}$ extends over  the whole $\moduli$ only as a sheaf and not a bundle.  

It is convenient to consider an arbitrary framing $\framing$ of $\E\in\moduli_0$; this is unique up to the $\P\GL_n$ action described above \eqref{quotPGL}. 
The framing $\framing$ of $\E$ allows us to regard a section of $\E^*\otimes \mathcal O(\scr D)$ as a vector valued meromorphic function $\psi$ on $\CC$ such that   
\be
\label{o12}
(\psi)\geq -\scr D,
\ee
and such that
\be
\label{o11}
\nf ^{-1}(z)\psi(z)\mbox{ is analytic at }\scr T,
\ee
where $\nf (z)$ is the polynomial normal form in $\DD$ as in the previous sections.

The following description of  connections follows essentially the ideas in \cite{kricheveriso}, and Theorem \ref{nomonothm} is also a reformulation, extended to the case of arbitrary degeneracy of the Tyurin data of Section 2 of loc. cit.
Denote by $\mathbb A_{(\E,\framing)}^{\scr D}$ the space of matrix-valued differentials ${\bf A}$ with poles only at $\scr T+\scr D$ and  such that
\be
\label{o13}
\nf ^{-1}{\bf A}\nf -\nf ^{-1}\d \nf \mbox{ is analytic at }\scr T,
\ee
\be
\label{o14}
{\bf A}+n_p\frac{\d z_p}{z_p}\mbox{ is analytic at }\scr D,
\ee
for any local coordinate $z_p$ near $p\in\scr D=\sum_p n_pp$, $z_p(p)=0$\footnote{If a point $p$ of $\scr D$ belongs to $\scr T$ as well, then the requirement for ${\bf A}$ is simply that
$$
\nf ^{-1}{\bf A}\nf -\nf ^{-1}\d \nf+n_p\frac{\d z_p}{z_p} \mbox{ is analytic at $p$}.
$$
}.
We have an identification $\conn_\E^{\scr D} {\simeq}{\mathbb  A}_{(\E,\framing)}^{\scr D}$, where the connection $\nabla \in \conn_\E^{\scr D}$ acts on the section $\psi$ in \eqref{o12}--\eqref{o11} as $ \nabla \psi = \d \psi- {\bf A}\psi$. The space $\mathbb A^{\scr D}_{(\E,\framing)}$ is an affine space over the space of Higgs fields $\Phi\in H^0({\rm End}(\E)\otimes\mathcal K)$; indeed if $\Phi$ is a matrix of differentials with $(\Phi)\geq -\scr T$ and $\nf ^{-1}\Phi \nf $ analytic at $\scr T$, then ${\bf A}+\Phi$ satisfies \eqref{o13}--\eqref{o14} whenever ${\bf A}$ does.

\begin{theorem}
\label{nomonothm}
For any ${\bf A}\in\mathbb A_{(\E,\framing)}^{\scr D}$ the differential equation $\d \Psi = {\bf A} \Psi$ for the matrix~$\Psi$ has apparent singularities at $\scr T$ and monodromy-free singularities at~$\mathscr D$. More precisely we have the following.
\begin{enumerate}
\item For any solution $\Psi$, in $\DD$ we have the local factorization $\Psi = \nf  H$ with $H$ analytic at $\scr T$.
\item For any $p\in \scr D$ of multiplicity $n_p$ and any solution $\Psi$, we have the local factorization $\Psi = z_p^{-n_p} H$, with $H$ locally analytic near $p$, for any local coordinate near $p$, $z_p(p)=0$. 
\end{enumerate}
\end{theorem}
{\bf Proof.}
For the first statement, define $H:=\nf ^{-1}\Psi$. Then $\d H={\bf A}_+ H$ where ${\bf A}_+:=\nf ^{-1}{\bf A}\nf  - \nf ^{-1} \d \nf $. By \eqref{o11} ${\bf A}_+$ is analytic at $\scr T$ and so $H$ is analytic at $\scr T$ as well. The proof of the second statement is completely analogous. 
\QED

In geometric terms, the holomorphic connection $\nabla\in\conn^{\scr D}_\E$, written $\nabla=\d-{\bf A}$ for an arbitrary framing $\framing$ and $\bf A\in\mathbb A^{\scr D}_{(\E,\framing)}$, equips the bundle $\E^*\otimes\O(\scr D)$ with a flat structure whose flat sections are the columns of a matrix solution of $\d\Psi={\bf A}\Psi$.
The monodromy of this equation computes the holonomy of the flat bundle, a fact which will be exploited below in Section~\ref{secPGLn}.

\paragraph{The Fay-differential.}
We now come to the first main contribution by constructing explicitly a holomorphically varying holomorphic connection on the complement of the non-abelian theta divisor.
For any $(\E,\framing)\in\wh\moduli_0$, let $\Cauchy_{\infty}$ be the associated Cauchy kernel introduced in Section~\ref{secCauchy}.
For $z,w$ local coordinates of two points in the same coordinate chart, define ${\bf F}(w)$ by the formula
\be
{\bf F}(w) := \lim _{z\to w} \le(\Cauchy_{\infty}(w,z)-\frac {\1 \d w}{w-z}\ri).
\label{defF}
\ee
One can check that if we write ${\bf F}(w) = F(w)\d w$,  then matrix $F(w)$ transforms as follows under a change of local coordinate $ \zeta= \zeta(w)$;
\be
\label{changeofvariableF}
 F(\zeta) =  \frac {\d w}{\d \zeta}  F(w)   +  \frac \1 2 \frac{\d}{\d \zeta} \ln \le( \frac {\d w}{\d \zeta} \ri).
\ee
Because of our initial assumptions on $\scr D$, there is a unique half-differential $h_{\scr D}$ (up to scalar) with divisor properties
\be
(h_{\scr D})\geq \scr D-\infty
\ee
and with ${\rm U}(1)$ multiplier system (see \cite{Fay73}). 

\begin{definition}
We define the {\bf Fay-differential} as 
\be
\label{defFDelta}
\F_{\scr D} ({p}) := \F({p}) - \1\d\ln h_\scr D({p}).
\ee
\end{definition}
From \eqref{changeofvariableF} it follows that $\F_{\scr D}$ is a matrix of differential forms with poles at $\scr T+\mathscr D$, while the pole at $\infty$ of $\F$ is canceled by the pole of $ \d \ln h_{\scr D}$; indeed, near $\infty$ the behaviour of $\Cauchy_\infty$ is given in a local coordinate $z$ with $z(\infty)=0$ by
\be
\frac{\Cauchy_\infty(w,z)}{\d w} = \frac {\1}{w-z} - \frac{\1}{w} + Q_\infty(w) + \mathrm O(z)
\ee
and hence $\F_{\scr D} (z) =  \mathrm{O}(1)\d z$ near $\infty$.

\begin{theorem}
\label{nosing}
{\bf [1]} The Fay-differential $\F_{\scr D}$ belongs to $\mathbb A_{(\E,\framing)}^{\scr D}$.

\noindent {\bf [2]} Under the $\P\GL_n$ action of change of framing we have $\F_{\scr D}\mapsto C\F_{\scr D}C^{-1}$. Therefore the Fay-differential provides a well-defined connection $\nabla_{\scr D}= \d-\F_{\scr D}\in\conn^{\scr D}_\E$.\\
\noindent
 {\bf [3]} Let $\wt{\scr D}$ be another nonspecial divisor of degree $g$. Then 
\be
\label{ganFaydiff}
\F_{\wt{\scr D}} = \F_{\scr D} + \Omega_{\scr D,\wt{\scr D}}\1,
\ee
where $ \Omega_{\scr D,\wt{\scr D}} = \d \ln \left({h_{\scr D}} / {h_{\wt{\scr D}}}\right)$ is the unique third kind differential
with purely imaginary periods and with simple poles at the points $p\in \scr D$  with residue $n_p$ and with simple poles at the points $q\in \wt{\scr D}$ with residue $-n_q$ (where $n_p, n_q$ are the multiplicities of the points in the corresponding divisors).
\end{theorem}

\noindent {\bf Proof.}
{\bf [1]}
By Corollary \ref{corollCauchy} the expression $ \nf ^{-1}(w) \Cauchy_\infty(w,z) \nf (z)$ is regular in $\DD$ except along the diagonal $w=z$. Here $w,z$ denote, by a slight abuse of notation, both the point and the coordinates. 
Near $w=z$ we have
\begin{align}
 \nf ^{-1}(w) \Cauchy_\infty(w,z) \nf (z)
&=\label{onehand} \frac {\1\d w}{w-z} - \nf ^{-1}(w) \nf '(w) \d w + \nf ^{-1} (w) \F(w) \nf (w) + {\mathrm O}(z-w).
\end{align}
On the other hand
\be
\label{otherhand}
\nf ^{-1}(w) \Cauchy_\infty(w,z) \nf (z) = \frac {\1\d w}{w-z} + \wt{\mathbf A}(w) + {\mathrm O}(z-w),
\ee
where $\wt{\mathbf A}(w)$ is analytic for $w\in\DD$, because $\nf ^{-1}(w) \Cauchy_\infty(w,z) \nf (z)$ has no singularities other than the pole on the diagonal.
By comparison of \eqref{onehand} and \eqref{otherhand} we deduce that
\begin{align}
- \nf ^{-1}(w) \nf '(w) \d w + \nf ^{-1} (w) \F(w) \nf (w)  = \wt{\mathbf A}(w)\mbox{ is analytic at }\scr T,
\label{www}
\end{align}
and so
\be
- \nf ^{-1}(w) \nf '(w) \d w + \nf ^{-1} (w) \F_{\scr D}(w) \nf (w)  = \wt{\mathbf A}(w)-\1\d\ln h_{\scr D}\mbox{ is analytic at }\scr T
\ee
and \eqref{o13} is established. 

To verify the condition  \eqref{o14} it suffices to note that the multiplicity $n_p$ of $p$ in $\scr D$ coincides with the order of vanishing of $h_{\scr D}$.

\noindent {\bf [2]} It follows from the easily established fact that the non-abelian Cauchy kernel transforms as $\Cauchy_\infty({q},{p})\mapsto C\Cauchy_\infty({q},{p})C^{-1}$ under the change of framing.

\noindent {\bf [3]} This statement follows immediately from the definition \eqref{defFDelta}. 
\QED

\begin{corollary}
\label{corollaryidentificationcotangconnection}
Identifying the fibers $T^*_{\E}\moduli_0\simeq H^0({\rm End}(\E)\otimes\mathcal K)$ of the bundle $T^*\moduli_0$ as in  Proposition~\ref{cotang}, the map
\begin{align}
\nn
\mathcal I:\ \conn^{\scr D}&\to T^*\moduli_0
\\
\label{identificationcotangconnection}
(\E,\nabla)&\mapsto(\E,\nabla_{\scr D}-\nabla)
\end{align}
is an isomorphism of holomorphic bundles over $\moduli_0$.
\end{corollary}

\subsection{Monodromy representation}
\label{secPGLn}

For any $(\E,\framing)\in\wh\moduli_0$ and any ${\bf A}\in\mathbb A_{(\E,\framing)}^{\scr D}$, a matrix of flat sections $\Psi$ is uniquely defined by
\be
\label{dpsi}
\d\Psi={\bf A}\Psi,\qquad \Psi(\infty)=\1.
\ee
According to Theorem~\ref{nomonothm} the  matrix of flat sections $\Psi$ extends to a meromorphic matrix function on the universal cover of $\CC$, with poles {\it only} at $\scr D$.  Under  analytic continuation along $\gamma\in \pi_1(\CC,\infty)$, $\Psi$ undergoes a transformation
\be
\label{eq:monodromymatrices}
\Psi(\gamma\cdot{p}) = \Psi({p}) M_\gamma^{-1}, \qquad M_\gamma\in\GL_n(\C).
\ee
Note that the homotopy class of $\gamma$ is in $\pi_1(\mathcal C)$ and not $\pi_1(\mathcal C\setminus \scr D)$ because $\Psi$ has trivial local monodromy at $\scr D$.
Therefore we obtain the 
map
\be
\label{framedmonodromy}
\wh{\mathcal M}_\E:\mathbb A^{\scr D}_{(\E,\framing)} \to {\rm Hom}( \pi_1(\CC,\infty), \GL_n(\C)).
\ee
The group $\P\GL_n(\C)$ acts on ${\rm Hom}\,(\pi_1(\CC,\infty),\GL_n(\C))$ via conjugation, and the quotient under this action is called the {\it character variety} which we denote by $\charvar$.

\begin{proposition}
The  map \eqref{framedmonodromy} for all $\E\in\moduli_0$ factors through a well defined map 
\be
\label{monodromy}
\mathcal M:\conn^{\scr D}\to \charvar,
\ee
independent of the framing of $\E$.
\end{proposition}

\noindent{\bf Proof.}
We note that under the change of framing the matrices $A,\Psi$ transform as follows
\be
{\bf A}\mapsto C{\bf A}C^{-1},\quad \Psi\mapsto C\Psi C^{-1},
\ee
the second of which enforces the normalization condition at $\infty$ in \eqref{dpsi}. 
This implies that the monodromy matrices transform as $M_\g\mapsto CM_\g C^{-1}$. 
Therefore the map $\mathbb A_{(\E,\framing)}^{\scr D}\to\charvar$ obtained by composing $\wh{\mathcal M}_\E$ in \eqref{framedmonodromy} with the quotient projection ${\rm Hom}\,(\pi_1(\CC,\infty),\GL_n(\C))\to\charvar$, factors through a well defined map $\mathcal M$ as in the statement.
\QED

Note that the monodromy map $\mathcal M$ in \eqref{monodromy} is holomorphic because the flat sections of $\d - \mathbf A$ (i.e. the solutions of \eqref{dpsi}) depend analytically on the coefficients due to the standard theorems of solutions of ODEs and these depend analytically on the Tyurin data.  This holomorphicity  is in contrast with the standard Narasimhan--Seshadri map \cite{Donaldson} which is only real-analytic (since it maps into the ${\rm U}_n$-character variety).

\br
\label{propPGL}
If $\scr D, \wt {\scr D}$ are  non-special  divisors of degree $g$, then the corresponding solutions $\Psi, \wt \Psi$ of the equation \eqref{dpsi}
are related  by, see \eqref{ganFaydiff},
\be
\wt \Psi ({p}) = c\frac{h_{\scr D}(p)} {h_{\wt{\scr D}}(p)}\Psi({p}),\qquad c=\lim_{p\to\infty}\frac{h_{\wt{\scr D}}(p)}{h_{\scr D}(p)}\in\C^*.
\ee
Hence the monodromy map $\mathcal M$ in \eqref{monodromy} depends on the choice of $\scr D$ only up to a scalar unitary representation of $\pi_1(\CC,\infty)$ depending only on the degree zero line bundle $\mathcal O(\scr D-\wt {\scr D})$; in particular $\mathcal M$ descends to a map to the $\P\GL_n$ character variety which is independent of the divisor $\scr D$.
\er

\section{Equivalence of the (complex) symplectic structure on \texorpdfstring{$T^*\moduli$}{T*V} and the Goldman bracket on \texorpdfstring{$\charvar$}{R}}\label{secEquiv}

In this section we prove our main result in Theorem~\ref{main}; it shows  that the pull-back by the monodromy map~$\mathcal M$ of the  Goldman symplectic structure on~$\charvar$  coincides with the canonical symplectic structure on $\conn^{\scr D}$. This latter symplectic structure is the one  induced by the identification of the $\conn^{\scr D}$ with the cotangent bundle $T^*\moduli_0$ via the map $\mathcal I$ in \eqref{identificationcotangconnection}.

We first observe that the non-degeneracy of the pairing \eqref{pairing1} implies that the Liouville one-form $\liouville$ on $T^*\moduli_0$ can be expressed as
\be
\label{locform}
\liouville(\E,\Phi) = \sum_{t\in \scr T} \res{t} \tr \left(\Phi \delta\nf\nf^{-1} \right) = \sum_{t\in \scr T} \res{t} \tr \left(\Phi \delta \Psi \Psi^{-1} \right),
\ee
where in the second equality we use the first part of Theorem~\ref{nomonothm} to write $\Psi=\nf H$ with $H$ analytic and analytically invertible in the neighbourhood of each point in  $\scr T$.
Here we denote by  $\delta$ the  differential with respect to the moduli of the bundle, as opposed to $\d$ representing the differential on $\CC$. The expression \eqref{locform} is independent of the framing.
Therefore, $\delta\liouville$ is the canonical (complex) symplectic structure on $T^*\moduli_0$.
 
The main tool for the proof of Theorem~\ref{main} is to rewrite the tautological pairing \eqref{pairing1} or equivalently \eqref{locform} as an integral over the canonical dissection of $\CC$. This rewriting yields the {\it Malgrange--Fay} one-form which we now define.

\paragraph{The Malgrange--Fay form.}
Consider  a canonical dissection of the Riemann surface $\CC$ along generators\\ $\a_1,\b_1,\dots,\a_g,\b_g$ of $\pi_1 (\mathcal C,\infty)$ satisfying
\be
\label{canrel}
 \alpha_1 \beta_1 \alpha_1^{-1} \beta_1^{-1}\cdots \alpha_g \beta_g \alpha_g^{-1} \beta_g^{-1}={\rm Id}.
\ee
 Let us denote by $\S$ the boundary of the fundamental polygon of this dissection of $\CC$, and by $\S_0$ the disjoint union of arcs obtained from $\S$ removing the points corresponding to $\infty$ (the vertices of the fundamental polygon).

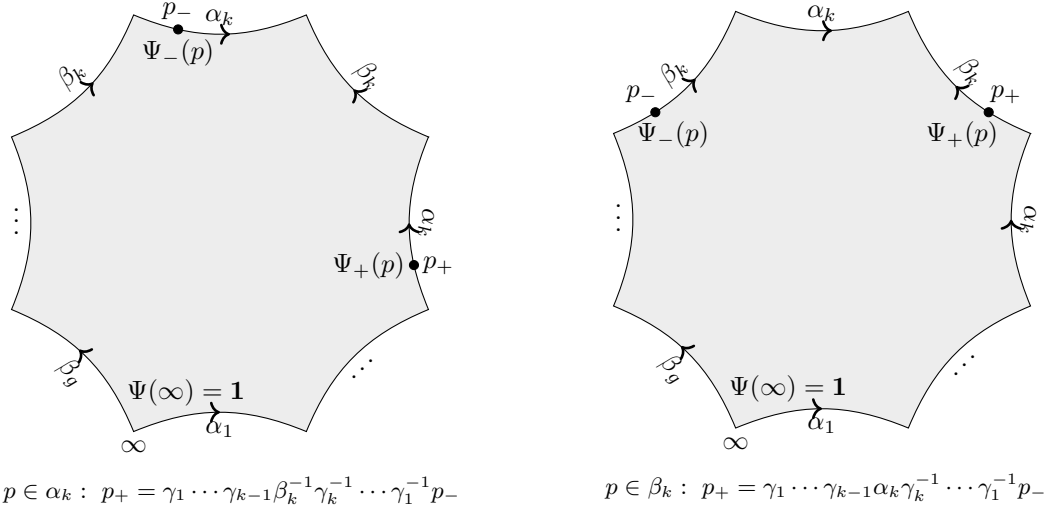
\begin{figure}[	t]
\centering
\hspace{1cm}
\begin{subfigure}{.45\textwidth}
\centering
\begin{tikzpicture}[scale=3]

\draw [fill=white!93!black 
, postaction={decorate,decoration={markings,mark=at position -2.5*0.125 with {\arrow[line width=1pt ]{>}}}}
, postaction={decorate,decoration={markings,mark=at position -.5*0.125 with {\arrow[line width=1pt ]{>}}}}
, postaction={decorate,decoration={markings,mark=at position .5*0.125 with {\arrow[line width=1pt ]{>}}}}
, postaction={decorate,decoration={markings,mark=at position 1.5*0.125 with {\arrow[line width=1pt ]{<}}}}
, postaction={decorate,decoration={markings,mark=at position 2.5*0.125 with {\arrow[line width=1pt ]{<}}}}
, postaction={decorate,decoration={markings,mark=at position 4.5*0.125 with {\arrow[line width=1pt ]{<}}}}
]
(0.5*45:1) to [in=-90+.5*45,out=135+.5*45] node[pos=0.5,sloped, above]{$\beta_k$}
      (1.5*45:1)   to [in=-45+.5*45,out=180+.5*45] node[pos=0.5,sloped, above]{$\a_k$}
      coordinate [pos=0.75] (meno)
      (2.5*45:1) to[in=.5*45,out=225+.5*45] node[pos=0.5,sloped, above]{$\b_k$}
      (3.5*45:1) to[in=45+.5*45,out=270+.5*45]node[pos=0.5,sloped, below]{$\cdots$}
      (4.5*45:1) to[in=90+.5*45,out=315+.5*45]node[pos=0.5,sloped, below]{$\b_g$}
      coordinate [pos=1] (infty)
      (5.5*45:1) to[in=135+.5*45,out=.5*45]node[pos=0.5,sloped, below]{$\a_1$}
      (6.5*45:1) to[in=180+.5*45,out=45+.5*45]node[pos=0.5,sloped, below]{$\cdots$}
      (7.5*45:1) to [in=225+.5*45,out=90+.5*45]node[pos=0.5,sloped, above]{$\a_k$}
      coordinate [pos=0.25] (piu)
      cycle;
\draw [fill](piu) circle[radius=0.02] node[right]{$p_+$} node[left]{$\Psi_+(p)$};
\draw [fill](meno) circle[radius=0.02] node[above]{$p_-$} node[below]{$\Psi_-(p)$};
\draw [fill](infty) circle[radius=0] node[below]{$\infty$};
\node at (-.15,-.75) {$\Psi(\infty)=\1$};
\end{tikzpicture}
\caption*{$p\in\a_k:\ p_+=\g_1\cdots\g_{k-1}\b_k^{-1}\g_k^{-1}\cdots\g_1^{-1}p_-$}
\end{subfigure}
\begin{subfigure}{.45\textwidth}
\centering
\begin{tikzpicture}[scale=3]

\draw [fill=white!93!black 
, postaction={decorate,decoration={markings,mark=at position -2.5*0.125 with {\arrow[line width=1pt ]{>}}}}
, postaction={decorate,decoration={markings,mark=at position -.5*0.125 with {\arrow[line width=1pt ]{>}}}}
, postaction={decorate,decoration={markings,mark=at position .5*0.125 with {\arrow[line width=1pt ]{>}}}}
, postaction={decorate,decoration={markings,mark=at position 1.5*0.125 with {\arrow[line width=1pt ]{<}}}}
, postaction={decorate,decoration={markings,mark=at position 2.5*0.125 with {\arrow[line width=1pt ]{<}}}}
, postaction={decorate,decoration={markings,mark=at position 4.5*0.125 with {\arrow[line width=1pt ]{<}}}}
]
(0.5*45:1) to [in=-90+.5*45,out=135+.5*45] node[pos=0.5,sloped, above]{$\beta_k$}
      coordinate [pos=0.25] (piu)
      (1.5*45:1) to [in=-45+.5*45,out=180+.5*45] node[pos=0.5,sloped, above]{$\a_k$}
      (2.5*45:1) to[in=.5*45,out=225+.5*45] node[pos=0.5,sloped,above]{$\b_k$}
      coordinate [pos=0.75] (meno)
      (3.5*45:1) to[in=45+.5*45,out=270+.5*45]node[pos=0.5,sloped, below]{$\cdots$}
      (4.5*45:1) to[in=90+.5*45,out=315+.5*45]node[pos=0.5,sloped, below]{$\b_g$}
      coordinate [pos=1] (infty)
      (5.5*45:1) to[in=135+.5*45,out=.5*45]node[pos=0.5,sloped, below]{$\a_1$}
      (6.5*45:1) to[in=180+.5*45,out=45+.5*45]node[pos=0.5,sloped, below]{$\cdots$}
      (7.5*45:1) to [in=225+.5*45,out=90+.5*45]node[pos=0.5,sloped, above]{$\a_k$}
      cycle;
\draw [fill](piu) circle[radius=0.02] node[above]{$\ \ \ \ p_+$} node[below]{$\Psi_+(p)\ \ \ \ \ \ $};
\draw [fill](meno) circle[radius=0.02] node[above]{$p_-\ \ \ $} node[below]{$\ \ \ \ \Psi_-(p)$};
\draw [fill](infty) circle[radius=0] node[below]{$\infty$};
\node at (-.15,-.75) {$\Psi(\infty)=\1$};
\end{tikzpicture}
\caption*{$p\in\b_k:\ p_+=\g_1\cdots\g_{k-1}\a_k\g_k^{-1}\cdots\g_1^{-1}p_-$}
\end{subfigure}
\caption{The fundamental flat section $\Psi$ satisfying \eqref{dpsi} is analytic in the interior of the fundamental polygon. Points on the boundaries are identified in pairs $p_\pm$, related by the indicated deck transformations ($\g_k:=\a_k\b_k\a_k^{-1}\b_k^{-1}$). The corresponding boundary values $\Psi_\pm(p)$ are then related by the associated monodromy transformation.}
\label{figjumps}
\end{figure}

Given a connection ${\bf A}\in \mathbb A^{\scr D}_{(\E,\framing)}$, by Theorem \ref{nomonothm} we shall regard a fundamental solution $\nabla\Psi=0$, see \eqref{dpsi}, as a single valued meromorphic function on $\CC\setminus\S$ with poles only at $\scr D$ and such that 
\begin{align}
\label{Psipoles}
z_p^{n_p} \Psi({p}) \mbox{ analytic at }p,\mbox{ for all }p \in \scr D,
\end{align}
where $n_p$ is the multiplicity of $p$ in $\scr D$ and $z_p$ a local coordinate at $p$, $z_p(p)=0$, compare with \eqref{o12}.
Note that since $\infty$ is the base-point of the dissection the  normalization $\Psi(\infty)=\1$ is intended in the sense that $\infty$ is approached within the sector bounded by $\b_g,\a_1$ (see Fig. \ref{figjumps}).

We denote the boundary values of $\Psi(p)$ at $\S_0$ as $\Psi_\pm(p)$, where $_+$ ($_-$, respectively)  denotes the boundary value from the left (right, respectively) of the oriented arcs $\a_k,\b_k$, see Fig. \ref{figjumps}.
Since points $p$ on $\S$ are identified in pairs, we have the following {\it jump relations}
\be
\label{RHPPsi} 
\Psi_+(p)=\Psi_-(p)J(p),\qquad p\in\S_0,
\ee
where the matrix $J$ is piecewise continuous on $\S_0$.
It is readily established that denoting $J_{\a_k},J_{\b_k}$ the values of $J(p)$ when $p\in\a_k,p\in\b_k$ (respectively) we have
\be
\label{jumpvsmonodromy}
J_{\a_k}=M_{\g_1\cdots\g_k\b_k\g_{k-1}^{-1}\cdots\g_1^{-1}},\qquad
J_{\b_k}=M_{\g_1\cdots \g_k\a_k^{-1} \g_{k-1}^{-1}\cdots \g_1^{-1}},
\ee
where $\g_k:=\a_k\b_k\a_k^{-1}\b_k^{-1}\in\pi_1(\CC,\infty)$, see Fig. \ref{figjumps}.

\begin{definition}
\label{defMalFay}
We define the {\bf Malgrange--Fay one-form} $\Malg$ on $\conn^{\scr D}$ by
\be
\label{MalgFay}
 \Malg(\E,\nabla):= \frac 1{2\pi\i}\int_{\Sigma}\tr\left(\Psi_-^{-1}\nabla_{\scr D}\Psi_- \delta J J^{-1}\right).
\ee
where $\nabla_\scr D=\d-\F_\scr D$ and $\Psi$ is the fundamental matrix of flat sections satisfying $\nabla\Psi=0,\Psi(\infty)=\1$, as in \eqref{dpsi}. Here $\S$ is the boundary of the canonical dissection, as above.
\end{definition}

The expression \eqref{MalgFay} is independent of the framing used to compute it.
Moreover, it does not depend on the choice of $\scr D$.
Indeed, if we change $\scr D $ to $\wt {\scr D}$, according to Remark \ref{propPGL} the matrix $\Psi$ is multiplied by a scalar factor that is independent of the moduli of $(\E,\framing)$.
The monodromy representation is also multiplied by a scalar  unitary representation of $\pi_1(\CC,\infty)$ (also independent of the moduli).
These factors cancel in the integrand of \eqref{MalgFay}.

The definition of the Malgrange--Fay form \eqref{MalgFay} comes partially from prior works \cite{BertoMalg,Bertocorr} with necessary adjustments to account for the nontrivial topology of the curve ${\mathcal C}$. 

\begin{proposition}
\label{propmalgfay}
The Malgrange--Fay form \eqref{MalgFay} is the pullback of the Liouville form $\liouville$ on $T^*\moduli_0$ along the identification \eqref{identificationcotangconnection} of Corollary \ref{corollaryidentificationcotangconnection};
\be
\Malg=\mathcal I^*\liouville.
\ee
\end{proposition}
\noindent {\bf Proof.}
From \eqref{RHPPsi} it follows that
\be
\label{jumppsi}
\Delta_{\Sigma} \le(\delta \Psi \Psi^{-1}\ri)  = \Psi _-\delta JJ^{-1} \Psi_-^{-1}.
\ee
where we denote by $\Delta_\Sigma$ the jump across $\S$,
\be
\label{jumpdef}
\Delta_\Sigma X({p}) = X_+({p})- X_-({p}),\ \ \ {p}\in \Sigma.
\ee
Therefore we rewrite \eqref{locform} as
\be
\label{rewriting}
\liouville = \sum_{t\in\scr T} \res{t} \tr \left(\Phi \delta \Psi \Psi^{-1} \right)
=\frac 1{2\pi\i}  \int_{\Sigma} \tr \left(\Phi\Delta_{\Sigma} (\delta \Psi \Psi^{-1} )\right)
=\frac 1{2\pi\i}  \int_{\Sigma} \tr \left(\Phi\Psi_- \delta J J^{-1}\Psi_-^{-1}\right).
\ee
In the second equality we have used Cauchy residue theorem, and noted that the poles of $\Psi$ at $\scr D$ do not contribute because $\delta \Psi \Psi^{-1}$ is holomorphic at $\scr D$ thanks to \eqref{Psipoles}, as the points of the divisor $\scr D$ do not enter the differential $\delta$.

The map $\mathcal I$ sends the connection $\nabla$ to the Higgs field $\Phi=\nabla_{\scr D}-\nabla$; hence from the last expression of \eqref{rewriting} we get
\be
\label{firstpullback}
\mathcal I^*\lambda=\frac 1{2\pi\i}\int_{\Sigma}\tr\left(\Psi_-^{-1}(\nabla_{\scr D}-\nabla)\Psi_- \delta J J^{-1}\right)
=\frac 1{2\pi\i}\int_{\Sigma}\tr\left(\Psi_-^{-1}\nabla_{\scr D}\Psi_- \delta J J^{-1}\right)=\Malg.
\ee
Here we have used that $\Psi$ are flat sections for the connection $\nabla$ and thus $\nabla\Psi=0$.
\QED

We next proceed with the computation of $\delta\Malg$.
To this end we recall the following notion of ``two-form associated to a graph'' from~\cite[Section~2]{BKGoldman}.
Let $\Sigma\subset \mathcal C$ be an arbitrary (locally finite)  graph consisting of smooth arcs meeting transversally at the vertices. We denote by ${\bf V}(\Sigma)$ the set of vertices of $\Sigma$, by ${\bf E}(\Sigma)$ the set of oriented edges of $\Sigma$ (namely, each edge appears twice in ${\bf E}(\S)$  with both possible orientations). 
An {\it admissible jump matrix on $\Sigma$} is a map $J: {\bf E}(\Sigma) \to \GL_n(\C)$ satisfying the following two requirements;
\begin{enumerate}
\item For each oriented edge $e\in {\bf E}(\Sigma)$,  we have $J(-e) = J(e)^{-1}$, where $-e$ denotes the oppositely oriented edge. 
\item For each vertex $v\in {\bf V}(\Sigma)$, consider the oriented edges which are incident to $v$ and oriented outwards; we enumerate them cyclically $e_1,\dots, e_{n_v}$ counterclockwise, where $n_v$ is the valence of $v$.
Then the requirement on $J$ (independent of the cyclical order chosen) is that 
\be
\label{nomono}
J(e_1)\cdots J(e_{n_v}) =\1.
\ee
\end{enumerate}

\begin{definition}[Two form $\Omega(\Sigma)$]
\label{defOmega}
Given a pair $(\Sigma, J)$ of an oriented graph and admissible jump matrix, we associate the following two-form 
\be
\label{twoformgraph}
\Omega(\Sigma):= \sum_{v\in {\bf V}(\Sigma)} \sum_{\ell=1}^{n_v}  \tr \left(
 K_\ell^{-1} \delta K_\ell \wedge J_\ell^{-1} \delta J_\ell 
\right),\qquad 
K_\ell := J_1\cdots J_{\ell},
\ee
where $J_\ell = J(e_\ell)$ are the jump matrices for the edges $e_{1},\dots, e_{n_v}$ incident at $v$, enumerated  counterclockwise starting from an arbitrary one and oriented away from $v$. 
\end{definition}

The form is invariant under cyclic permutations of the labels of the incident edges at a vertex, thanks to \eqref{nomono}.
It is shown in \cite{BKGoldman} that $\Omega$ is invariant under natural transformations of the graph like edge contractions etc. We will not need them specifically here and thus do not recall the details. 
 
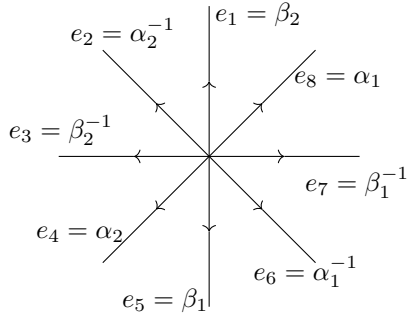
\begin{figure}[t]
\centering
\hspace{1cm}
\begin{tikzpicture}[scale=1]
\begin{scope}[rotate=45]
\draw[->] (0,0) to (0:1);
\draw[->] (0,0) to (45:1);
\draw[->] (0,0) to (90:1);
\draw[->] (0,0) to (135:1);
\draw[->] (0,0) to (180:1);
\draw[->] (0,0) to (-135:1);
\draw[->] (0,0) to (-90:1);
\draw[->] (0,0) to (-45:1);
\draw[-] (0:1) to (0:2);
\draw[-] (45:1) to (45:2);
\draw[-] (90:1) to (90:2);
\draw[-] (135:1) to (135:2);
\draw[-] (180:1) to (180:2);
\draw[-] (-135:1) to (-135:2);
\draw[-] (-90:1) to (-90:2);
\draw[-] (-45:1) to (-45:2);
\node at (-14:2) {$e_8=\a_1$};
\node at (-55:2) {$e_7=\b_1^{-1}$};
\node at (-95:2) {$e_6=\a_1^{-1}$};
\node at (-152:2) {$e_5=\b_1$};
\node at (-194:2) {$e_4=\a_2$};
\node at (-235:2) {$e_3=\b_2^{-1}$};
\node at (-280:2) {$e_2=\a_2^{-1}$};
\node at (-334:2) {$e_1=\b_2$};
\end{scope}
\end{tikzpicture}
\caption{Half edges of $\S$ incident at $v=\infty$ for $g=2$.}
\label{fighalfedges}
\end{figure}

In our case we take $\S$ to be the canonical dissection. Then the jump matrix $J$ in \eqref{RHPPsi}, given in \eqref{jumpvsmonodromy}, are an admissible jump on $\S$. Indeed the graph $\S$ has only one vertex at $\infty$ and so we label the half edges $e_1,\dots,e_{4g}$ of $\S$, in counterclockwise sense, so that $e_{1}$ corresponds to $\b_g$, $e_2$ to $\a_g^{-1}$, $e_3$ to $\b_g^{-1}$, $e_4$ to $\a_g$, and so on as in Fig. \ref{fighalfedges}; denoting $J_\ell=J(e_\ell)$ we have
\be
J_{4(g-\ell)+1}=J_{\b_\ell},\quad J_{4(g-\ell)+2}=J_{\a_\ell}^{-1},\quad J_{4(g-\ell)+3}=J_{\b_\ell}^{-1},\quad J_{4(g-\ell)+4}=J_{\a_\ell},\qquad \ell=1,\dots,g.
\label{714j}
\ee

\begin{shaded}
\begin{theorem}
\label{thmmalgfay}
The exterior derivative of the Malgrange--Fay form $\Malg$ of Def. \ref{defMalFay} is given in terms of the monodromy matrices by 
\be
\label{bbbb}
-{4\pi\i}\, \delta \Malg = \Omega(\Sigma)  =\sum_{\ell=1}^{4g} \tr \left(K_\ell^{-1} \delta K_\ell \wedge J_\ell^{-1} \delta J_\ell\right),\qquad 
K_\ell := J_1\cdots J_{\ell},
\ee
where $\Sigma$ is the canonical dissection and $\Omega(\S)$ is in \eqref{twoformgraph}. The matrices $J$ are given by the holonomy of the flat connection $\nabla$, see \eqref{jumpvsmonodromy} and \eqref{714j}.
\end{theorem}
\end{shaded}

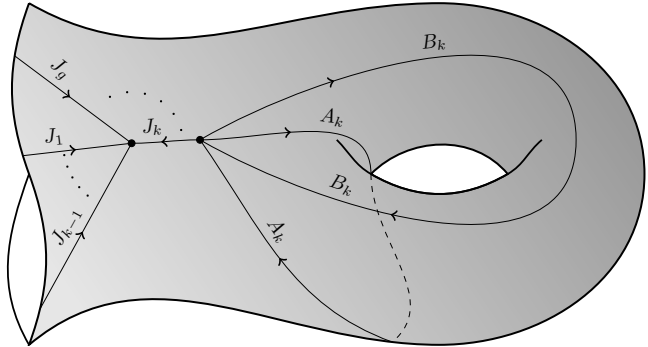
\begin{figure}[t]
\begin{center}
\begin{tikzpicture}[scale=2, x=0.3cm, y=0.3cm, xshift=-2]
\begin{scope}[rotate=0]
\draw [fill] (0,0) circle[radius =0.03];

\coordinate  (p1) at (-5,4);
\coordinate (p2) at  (7,4);
\coordinate (p3) at (13,-1);
\coordinate (p4) at (7,-6);
\coordinate (p5) at (-5, -6);
\coordinate (p6) at (-5, -1);
\coordinate (q1) at (5,-1);
\coordinate (q11) at (4,0);
\coordinate (q2) at (9,-1);
\coordinate (q22) at (10,0);
\coordinate(s1) at (11, 0);

\coordinate (t1) at (-3,-1.3);
\coordinate (t2) at (-2,-3);
\coordinate (tn) at (-0,-3.3);
\coordinate(v0) at (-2,-0.1);

\draw [
line width=1,
top color= white!60!black, 
bottom color= white!90!black, 
shading=axis, shading angle=-70.02]
(p1) to [out=-30, in=180] (p2) 
to [out=0, in=90] (p3)
to [out=-90, in =0] (p4)
to [out=180, in=40, looseness=1] coordinate [ pos=0.1] (s2) (p5)
to [out=70, in=-70] 
coordinate[pos=0.2] (ex1)
coordinate[pos=0.6] (ex2)
(p6) 
to [out = 110, in=-110]  
coordinate[pos=0.1] (ex3)
coordinate[pos=0.6] (ex4)
coordinate[pos=0.7] (ex5)
coordinate[pos=0.8] (ex6)
(p1);

\draw [line width=0.9] (p5) to [out=120, in=-110] (p6);

\draw [fill=white, line width=1]
(q1) to [out=45, in=130] (q2) to [out= 210, in=-30] (q1);

\draw [line width =0.9] 
(q11) to [out=-45, in=150] (q1) to [out=-30, in=-150] (q2) to [out=30, in = -140] (q22);

\draw [fill] (0,0) circle[radius=0.1];
\draw [fill] (v0) circle[radius=0.1];
\draw 
[    postaction={decorate,decoration={markings,mark=at position 0.1666 with {\arrow[line width=1pt ]{>}}}},
postaction={decorate,decoration={markings,mark=at position 0.7666 with {\arrow[line width=1pt ]{>}}}},
]
(0,0) to [out=30, in =90] 
node [pos=0.5, above, sloped] {$B_k$}
(s1) to [out=-90, in =-30] 
node [pos=0.7, above, sloped] {$B_k$}
(0,0);

\draw [postaction={decorate,decoration={markings,mark=at position 0.4666 with {\arrow[line width=1pt ]{>}}}}] (0,0) to [out=0, in=90] node [pos=0.6, above, sloped] {$A_k$} (q1);
\draw [dashed] (q1) to [out=-90, in = 40] (s2);
\draw [postaction={decorate,decoration={markings,mark=at position 0.5 with {\arrow[line width=1pt ]{>}}}}]
(s2) to [out=160, in=-55] node [pos=0.6, above, sloped] {$A_k$} (0,0);

\foreach \vv in {190,200,210,220,230} {
\draw ($(v0)+(\vv:2)$) circle [radius=0.15pt];
}

\foreach \vv in {15,35,55,75,95,115} {
\draw ($(v0)+(\vv:1.5)$) circle [radius=0.15pt];
}

\draw  [postaction={decorate,decoration={markings,mark=at position 0.566 with {\arrow[line width=1pt ]{>}}}}]   (0,0) to [out =180 , in =0] node [pos=0.7, above, sloped] {$J_k$} (v0) ;

\draw  [postaction={decorate,decoration={markings,mark=at position 0.566 with {\arrow[line width=1pt ]{<}}}}]   (v0) to 
 node [pos=0.55, above, sloped] {$J_{k-1}$} (ex1) ;
\draw  [postaction={decorate,decoration={markings,mark=at position 0.566 with {\arrow[line width=1pt ]{<}}}}]   (v0) to 
 node [pos=0.7, above, sloped] {$J_1$} (ex3) ;

\draw  [postaction={decorate,decoration={markings,mark=at position 0.566 with {\arrow[line width=1pt ]{<}}}}]   (v0) to 
 node [pos=0.7, above, sloped] {$J_g$} (ex5) ;
\end{scope}
\end{tikzpicture}
\end{center}
\caption{The Krichever graph $\Sigma_{\text{Kr}}$ with indicated the jump matrices as in \cite{kricheveriso} for comparison, so that $J_k = B_k^{-1} A_{k}^{-1} B_k A_k$.  The two-form associated to this graph following Def. \ref{defOmega} is equivalent, after graph contractions, to the form in \cite{AlMal} and hence to the Goldman symplectic form. The matrices that Krichever denotes by $A_k, B_k$  coincide with what we denote here by {$J_{\alpha_k}, J_{\beta_k}$}.}
\label{figkrich}
\end{figure}

This theorem is central to our main result on the Goldman symplectic structure, and its proof is deferred to  Appendix \ref{proofmain}. 
{
\begin{remark}
\label{remkri}
It is appropriate to remark here the connection with the results  in \cite{krichever,kricheveriso}, where the author considers the ``universal symplectic structure": this is indeed the same, with the proper identifications, as the exterior differential of $\liouville$ in \eqref{locform}.  In the same papers the author shows that under the (extended)  monodromy map, the universal two-form is mapped to an expression (Theorem 6.1 in loc.cit.) which is, {\it prima-facie} different from \eqref{bbbb}.

We now indicate  how to show that it is equal (up to sign) to \eqref{bbbb}, which is in principle a  non-trivial verification because the two expressions match only on the manifold of constraints given by the fundamental relation in the homotopy group. From (6.9) of \cite{kricheveriso} it is possible to ``reverse engineer'' the graph $\Sigma_{\text{Kr}}$ and the jump matrices $J$ that yield Krichever's formula (6.9). This is shown in Fig. \ref{figkrich}.
Then  the expression (6.9) in loc. cit. \footnote{This is the case without additional punctures.} can be recognized to be equal (after a bit of simplification)   to the two form $\Omega(\Sigma_{\tiny \text{\tiny Kr}})$ of Def. \ref{defOmega} for the graph $\Sigma_{\tiny \text{\tiny Kr}}$ indicated in Fig. \ref{figkrich}. The fact that this coincides with the prior result of Alekseev and Malkin \cite{AlMal}, and hence with the Goldman symplectic form, is then a consequence of the invariance of the two form $\Omega(\Sigma)$ under standard transformations like edge-contractions, proved in \cite{BKGoldman}. In fact the graph $\Sigma$ of the canonical dissection  is the contraction of all the edges carrying the matrices $J_k = B_k^{-1}A_k^{-1} B_k A_k$ of the graph in Fig. \ref{figkrich}.
\end{remark}}

\paragraph{Equivalence with the Goldman bracket.}
The (complex) Goldman bracket \cite{Goldman} is a Poisson bracket $\{,\}_G$ on the character variety $\charvar$; for two arbitrary $\gamma, \wt \gamma\in \pi_1(\mathcal C,\infty)$ it is defined by
\be
\label{Goldman}
\bigg\{\tr M_\gamma , \tr M_{\wt \gamma} \bigg\}_G = \sum_{p\in \gamma\cap\wt \gamma} \nu(p) \tr (M_{\gamma_p\wt \gamma})
\ee 
where $\gamma,\wt\gamma$ are representative of the corresponding homotopy classes which meet transversally, so that $\g\cap\wt\g$ is a finite set, $\nu(p)=\pm 1$ is the intersection index of $\gamma$ and $\wt \gamma$ at $p$, and $\g_p\wt\g$ is the homotopy class of the path obtained traversing $\g$ starting from  the intersection point $p$, then traversing the loop $\wt\g$.

In the case under consideration, where $\mathcal C$ is without punctures and without boundary, the Poisson bracket \eqref{Goldman} is known to be non-degenerate; the corresponding symplectic form $\varpi_{G}$ can be defined by the relation
$
\{f,g\}_G=\varpi_G(X_f,X_g),
$
denoting $X_f=\{,f\}_G$ the Hamiltonian vector fields, for any smooth functions $f,g$ on $\charvar$.

\begin{shaded}
\begin{theorem}
\label{main} 
The pull-back to $\conn^{\scr D}$ of the Goldman symplectic structure $\varpi_G$ along the monodromy map $\mathcal M:\conn^{\scr D} \to\charvar$  is proportional to  the canonical symplectic structure on $\conn^{\scr D}$ induced by the isomorphism $\mathcal I:\conn^{\scr D}\to T^*\moduli_0$;
\be
-2\pi\i\,\mathcal I^* \delta \liouville  = \mathcal M^*\varpi_G.
\ee
\end{theorem}
\end{shaded}

\noindent {\bf Proof.}
The symplectic form, $\varpi_G$, associated to the Goldman Poisson structure was provided in \cite{AlMal}. In \cite{BKGoldman}, Theorem~3.1 it is shown that the expression of Alekseev-Malkin, after appropriate transformations, equals $\frac 1 2 \Omega(\Sigma)$ where $\Sigma$ is the canonical dissection, and $\Omega(\S)$ is the two-form of Definition \ref{defOmega}; namely we have
\be
\frac 12\Omega(\S)=\mathcal M^*\varpi_G.
\ee
The proof is completed by the chain of equalities
\begin{multline*}
\shoveleft {
-2\pi\i\,\mathcal I^*\delta\liouville\mathop{=}^{
\hbox{
\resizebox{20pt}{!}{{  Prop.~\ref{propmalgfay}} }
}
}
 -2\pi\i\,\delta\Malg
\mathop{=}^{\hbox{
\resizebox{20pt}{!}{{  Thm. \ref{thmmalgfay}}}
}}
\frac 1 2 \Omega(\Sigma)=\mathcal M^*\varpi_G .
\hspace{3cm} \qquad \qquad \blacksquare}
\end{multline*}

\begin{remark}[Relationship with Fay's \cite{Fay92}]
\label{compaFay}

In \cite{Fay92} Fay introduced, in a slightly different context, a similar one-form. The main difference is that Fay works with flat vector bundles tensored by a root of the canonical bundle, and hence he works on vector bundles of degree $n(g-1)$ so that Riemann--Roch's theorem states $h^0=h^1$ and both are generically zero.
\end{remark}

\section{Singularities of the Malgrange--Fay form  along the non-abelian theta divisor}

The Cauchy kernel is ill-defined on the non-abelian theta divisor $\wh{\thetadiv}$, and hence we do not have a way to extend the reference connection $\F_\scr D$ over $\wh \thetadiv$.
In this last section we show that the Malgrange--Fay form $\Malg$ on the bundle of connections $\conn^{\scr D}$ has a first order pole with residue given  by the index $h^1(\E)$.

This is the content of the following theorem, which is an analogue of Theorems 3 and 4 in \cite{Fay92}.

More precisely, we consider a holomorphic family $\E_{\epsilon}$, for $|\epsilon|<1$, such that $\E_{0}\in\thetadiv$. 
On the family  $\E_{\epsilon}$ we consider  framings $\framing_{\epsilon}$, and the associated polynomial normal forms $\nf_{\epsilon}$, Tyurin divisors $\mathscr T_{\epsilon}$, and Brill--Noether--Tyurin matrices $\mathbb T_{\epsilon}$. 

 Recalling that $\thetadiv$ is defined by $\det \mathbb T=0$, we call the family  $\E_{\epsilon}$ {\it transversal} to $\thetadiv$ if  
\be\label{transversality}
\res{\epsilon=0}\d\ln\det\mathbb T_{\epsilon}=h^1(\E_{0}).
\ee

We correspondingly take any family of connections $\nabla_\epsilon= \d - {\bf A}_\epsilon$ such that ${\bf A}_\epsilon$  satisfies the conditions \eqref{o13}, \eqref{o14} identically in $\epsilon$ and is  holomorphic for  $|\epsilon|<1$. 

\begin{theorem}
\label{theoremFay}
Denote by $\Malg_\epsilon$ the pull-back of $\Malg$ along a transversal holomorphic  family $(\E_{\epsilon},\nabla_\epsilon) \in \conn^{\scr D}$, $0<|\epsilon|<1$. Then $\Malg_\epsilon$ has a simple pole at $\epsilon=0$, with residue $h^1(\E_0)$, namely 
\be
\Malg_\epsilon=h^1(\E_0)\frac {\d \epsilon}\epsilon + \mathcal O(1).
\ee
\end{theorem}

Before proving the theorem we need some remarks and a lemma.
Equation 
\eqref{transversality} holds true if and only if $\dot{\mathbb T}_0:=\left.\pa_\epsilon\mathbb T_{\epsilon}\right|_{\epsilon=0}$ restricts to a non-degenerate pairing between left and right null-spaces of $\mathbb T_{0}$. Thanks to the identifications of the null-spaces of the Brill--Noether--Tyurin matrix provided in Proposition~\ref{propdual}, such restriction of $\dot{\mathbb T}_{0}$ can be expressed (up to an inessential sign) by the pairing
\be
\label{pairingres}
H^0(\E_0(-\infty))\otimes H^0(\E_0^*\otimes\mathcal K)\to \C:(\bs r,\bs\eta)\mapsto\left\langle\bs r,\bs\eta\right\rangle_0
=\sum_{t\in\scr T_0}\res{t}\bs r^t \dot \nf _0 \nf ^{-1}_0 \bs\eta,
\ee
where $\dot \nf _0=\left.\pa_\epsilon \nf _\epsilon\right|_{\epsilon=0}$.

\begin{lemma}
\label{lemmatransversality}
{\bf [1]}For any $\E_0\in\thetadiv$ and an arbitrary framing $\framing_0$, there exists a transversal deformation $\left(\E_{\epsilon},\framing_{\epsilon}\right)$, i.e. a deformation such that the pairing \eqref{pairingres} is non-degenerate. 

\noindent{\bf [2]}
Let $\left(\E_{\epsilon},\framing_{\epsilon}\right)$ be a transversal deformation of $(\E_0,\framing_0)\in\wh\thetadiv$ and let us denote $\Cauchy_\infty(q,p;\epsilon)$ the Cauchy kernel associated to $\left(\E_{\epsilon},\framing_{\epsilon}\right)$. Let $k=h^1(\E_0)\geq 1$ and ${\bf r}_1,\dots,{\bf r}_k$ and $\bs\eta_1,\dots,\bs\eta_k$ be arbitrary bases of $H^0(\E_0(-\infty))$ and $H^0(\E^*_0\otimes\mathcal K)$.
Then $\Cauchy_\infty(q,p;\epsilon)$ has a simple pole at $\epsilon=0$, whose residue is the kernel of rank $k$ given by
\be
\label{limK}
\lim_{\epsilon\to 0}\,\epsilon\, \Cauchy_\infty({p},{q};\epsilon) =
-
 \sum_{a,b=1}^k Q^{ab} {\bs \eta}_{a}({p}){\bf r}_{b}^t({q})
\ee 
where $Q^{ab}$ are the entries of the inverse $Q^{-1}$ of the Gram matrix $Q=(Q_{ab})_{a,b=1}^k$, $Q_{ab}=\langle{\bf r}_a,\bs\eta_b\rangle_0$.
\end{lemma}

\noindent{\bf Proof.}
\noindent{\bf [1]} The proof follows along the lines of an argument by Fay \cite{Fay92} adapted to our setting. We start by observing that if ${\bf r}\in H^0(\E_0(-\infty))$ and $\bs\eta\in H^0(\E_0^*\otimes\mathcal K)$ then $\Phi={\bs\eta}{\bf r}^t$ is a Higgs field, see \eqref{boundaryhiggs}.
The pairing \eqref{pairing1} with $\Phi={\bs\eta}{\bf r}^t$ coincides with \eqref{pairingres};
\be
(\Phi,\dot \nf _0)=\sum_{t \in \scr T_0} \res{t} \tr \left(\Phi \dot \nf _0\nf ^{-1}_0\right)=
\sum_{t \in \scr T} \res{t} {\bf r}^t \dot \nf _0\nf ^{-1}_0 \bs\eta
=\left\langle{\bf r},\bs\eta\right\rangle_0.
\ee
For any $[{\bf r}]\in\mathbb P H^0(\E_0(-\infty))$ let $V_{[{\bf r}]}$ be the subspace of $T_{\E_0}\wh\moduli$ defined by $\langle{\bf r},\cdot\rangle_0=0$; equivalently, denoting $k=h^1(\E_0)\geq 1$ and fixing a basis $\bs\eta_1,\dots,\bs\eta_k$ of $H^0(\E_0\otimes\mathcal K)$, then the space $V_{[{\bf r}]}$ is defined by the $k$ linear equations $({\bs\eta}_i{\bf r}^t,\dot \nf _0)=0$ for $i=1,\dots,k$.
Since $(,)$ is a nondegenerate pairing by Proposition \ref{propnondegenerate}, these $k$ linear equations are independent and we have $\dim V_{[{\bf r}]}=n^2g-k$.
Therefore the set of deformations $\dot \nf _0$ for which $\langle,\rangle_0$ is degenerate coincides with $\bigcup_{[{\bf r}]\in\mathbb P H^0(\E_0(-\infty))}V_{[{\bf r}]}$, which is of dimension at most $(n^2g-k)+(k-1)$ and thus cannot cover $T_{\E_0}\wh\moduli$.
Hence there exists $\dot \nf _0$ such that $\langle,\rangle_0$ is non-degenerate.

\noindent{\bf [2]} Looking at \eqref{expressioncauchydet} it is clear that $\Cauchy_\infty(q,p;\epsilon)$ has at most a simple pole at $\epsilon=0$; indeed due to the transversality condition $\res{\epsilon=0}\d\ln\det\mathbb{T}_{\epsilon}=k$ it is easy to conclude that the numerator in \eqref{expressioncauchydet} vanishes at $\epsilon=0$ of order at least $k-1$, while the denominator vanishes of order exactly $k$.
Hence we can write $\Cauchy_\infty({q},{p};\epsilon)$ as a Laurent expansion:
\be
\label{Laurent}
\Cauchy_\infty({q},{p};\epsilon)
=\sum_{j\geq -1}\Cauchy_\infty^{[j]}(q,p)\epsilon^j.
\ee
Since the Cauchy kernel has a simple pole along the diagonal with residue $\1$ (independent of $\epsilon$) it follows that the kernel $\Cauchy_\infty^{[0]}({q},{p})$ has a simple pole with residue $\1$ along ${q}={p}$ while  $\Cauchy_\infty^{[j]}$, $j\neq 0$ are analytic at ${q}={p}$. 
As a consequence of Corollary \ref{corollCauchy} we  also have that:
\begin{itemize}
\item the rows of $\Cauchy_\infty^{[-1]}(q,p)$, as meromorphic vector functions of $p$, are (transposition of) holomorphic sections of $\E_0(-\infty)$;
\item the columns of $\Cauchy_\infty^{[-1]}(q,p)$, as differentials in $q$, are holomorphic sections of $\E_0^*\otimes\mathcal K$.
\end{itemize}
Therefore we have
\be
\Cauchy_\infty^{[-1]}(q,p)=-\sum_{a,b=1}^k Q^{ab}\bs\eta_{a}(q)\mathbf{r}^t_{b}(p)
\ee
for some matrix $(Q^{ab})$ (the minus sign is for convenience).
To complete the proof we have to show that $Q^{ab}$ are the entries of the inverse of the matrix $Q_{ab}=\langle{\bf r}_a,\bs\eta_b\rangle_0$. To this end we first note that
\be
\label{uniformconvergence}
\Cauchy_\infty^{[-1]}(q,p)\dot \nf _0(p)=\lim_{\epsilon\to 0}\Cauchy_\infty(q,p;\epsilon)\left(\nf (p;\epsilon)-\nf (p;0)\right).
\ee
Choose now  $\bs\eta\in H^0(\E_0^*\otimes\mathcal K)$ and  assume $q\in\CC\setminus\DD$:  applying the Cauchy theorem we then obtain
\begin{align}
\nonumber
\left\langle\Cauchy_\infty^{[-1]}(q,\cdot),\bs\eta\right\rangle_0
&=\sum_{t\in\scr T_0}\res{p=t}\Cauchy_\infty^{[-1]}(q,p)\dot \nf _0(p)\nf _0^{-1}(p)\bs\eta(p)
=\frac 1{2\pi\i}\oint_{p\in\partial\DD}\Cauchy_\infty^{[-1]}(q,p)\dot \nf _0(p)\nf _0^{-1}(p)\bs\eta(p)
\\
&=\frac 1{2\pi\i}\lim_{\epsilon\to 0}\left(\oint_{p\in\partial\DD}
\Cauchy_\infty(q,p;\epsilon)\nf _\epsilon(p)\nf _0^{-1}(p)\bs\eta(p)-\Cauchy_\infty(q,p;\epsilon)\bs\eta(p)\right),
\end{align}
where exchanging the limit and the integral is legitimate because \eqref{uniformconvergence} is uniform with respect to  $p\in \pa \DD$.
The expression $\Cauchy_\infty(q,p;\epsilon)\nf (p;\epsilon)\nf _0^{-1}(p)\bs\eta(p)$ is regular for $p\in\DD$ (as $\nf _0^{-1}(p)\bs\eta(p)$ is regular for $p\to\scr T_0$) and so does not contribute to the integral. The remaining part is
\be
-\frac 1{2\pi\i}\oint_{p\in\partial\DD}\Cauchy_\infty(q,p;\epsilon)\bs\eta(p)=\res{p=q}\Cauchy_\infty(q,p;\epsilon)\bs\eta(p)=-\bs\eta(q).
\ee
Therefore we have shown that for all $q\in\CC\setminus\DD$ we have $\left\langle\Cauchy_\infty^{[-1]}(q,\cdot),\bs\eta\right\rangle_0=-\bs\eta(q)$. We now apply this identity to the basis $\{\bs\eta_c, \ c=1,\dots, k\}$ of $H^0(\E^*\otimes \mathcal K)$;
\be
\left\langle\Cauchy_\infty^{[-1]}(q,\cdot),\bs\eta_c\right\rangle_0=
-
\sum_{a,b=1}^k Q^{ab}\bs\eta_{a}(q)\left\langle{\bf r}_b,\bs\eta_c\right\rangle_0=
-
\sum_{a,b=1}^k Q^{ab}Q_{bc}\bs\eta_{a}(q)=-\bs\eta_{c}(q),
\ee
yielding $\sum_{b=1}^k Q^{ab}Q_{bc}=\delta^a_c$, and the proof is complete.
\QED

We now are in a position to prove Theorem \ref{theoremFay}.
\\

\noindent {\bf Proof of Theorem \ref{theoremFay}.}
Using Proposition \ref{propmalgfay} and formula \eqref{locform} can rewrite the pullback of the Malgrange form as
\be
\Malg_\epsilon = \sum_{t\in \scr T} \res{t} \tr \left(\Phi_\epsilon \dot \nf_\epsilon \nf_\epsilon^{-1} \right)\d\epsilon,
\label{811}
\ee  
where $\nf _\epsilon$ are the polynomial normal forms associated to the framed bundles $(\E_\epsilon,\framing_\epsilon)$ as above, and the family of connections $\nabla_\epsilon$ are trivialized as $\d-\mathbf A_\epsilon$ with ${\mathbf A}_\epsilon=\Phi_\epsilon + \F_{\scr D,\epsilon}$ where $\F_{\scr D, \epsilon}$ denotes the Fay connection on $\E_\epsilon$ and has a pole at $\epsilon=0$ 
according to Lemma \ref{lemmatransversality}, and in particular formula \eqref{limK}. The singular part is  of the form 
\be
\F_{\scr D}(z;\epsilon) = - \frac 1{\epsilon}\sum_{a,b} Q^{ab} {\bs \eta}_{a}(z) {\bf r}_{b}^t(z)  + {\rm O}(1) =: \frac 1 \epsilon \Phi^{[-1]}(z) + {\rm O}(1),
\ee
where ${\rm O}(1)$ denotes terms regular at $\epsilon=0$.
In particular the coefficient $\Phi^{[-1]}$ is a certain Higgs field on $\E_0$.
Since, by construction, $\mathbf A_\epsilon$ is chosen holomorphic in $\epsilon$, near $\epsilon=0$ we have
\be
\Phi_\epsilon=-\frac 1\epsilon\Phi^{[-1]}(z)+{\rm O}(1).
\ee
Inserting this expression into \eqref{811} we obtain
\be
\Malg_\epsilon = -\frac {\d \epsilon}{\epsilon} \sum_{t\in \scr T} \res{t} \tr \left(\Phi^{[-1]} \dot\nf_0\nf_0^{-1} \right) + {\rm O}(1) = \frac {\d \epsilon}{\epsilon} \sum_{a,b=1}^{h^1(\E_0)} Q^{ab} \overbrace{\left\langle{\bf r}_b,\bs\eta_a\right\rangle_0}^{=Q_{ba}} + {\rm O}(1)  = h^1(\E_0)\frac { \d \epsilon}{\epsilon} + {\rm O}(1).
\ee
This concludes the proof.
\QED

\section{Comments and open problems}
In this paper we have extensively employed the holomorphic parametrization of vector bundles in terms of Tyurin data and the notion of non-abelian Cauchy kernel; to conclude, we comment on a few future directions and interesting open questions that arise from these constructions.
\begin{itemize}
\item[-] The connection $\F_\scr D$, as a function of the moduli, is meromorphic with poles along the non-abelian theta divisor. 
In the setting of projective connections described in \cite{BKN} it compares with the Wirtinger projective connection. 
 
On the other hand in \cite{BKN} the central object was the Bergman projective connection, which is holomorphic on the Teichmuller space but it is not single-valued on the moduli space $\mathcal M_{g}$. The question is then whether there is a similar reference connection  which is holomorphic in the moduli of the bundles also on the theta divisor but not necessarily single valued.
\item[-] 
 We would like to define the notion of tau function in a similar spirit as  \cite{Fay92};
 since  $\delta \Malg$ is an expression  in terms of the monodromy matrices alone, one can find a local antiderivative $\vartheta$ in terms of the same matrices, locally analytic at the theta divisor. Thus $\Malg - \vartheta$ is a locally defined closed one form which defines a ``tau'' function as $ \delta \ln \tau = \Malg-\vartheta$. By Theorem \ref{theoremFay} this function vanishes at the points $\E$ in the theta divisor $\thetadiv$ to order $h^1(\E)$.   
The concrete obstacle here is that one needs to construct an {\it explicit} local potential for the Goldman symplectic structure. 
 
\item[-] In perspective, one could consider next {\it meromorphic} connections;  this is, implicitly, the setting of \cite{kricheveriso}. However in loc. cit. the notion of tau function was not investigated, and thus neither its role as generating function for change of Lagrangian polarization.  The natural conjecture is now that one can introduce a tau function by combining the logic indicated in the previous points. The resulting tau function should have  the properties of vanishing on the non-abelian theta divisor.
\item[-] 
A natural question arises in the context of {\it Hecke modifications} (see e.g. \cite{EnRub}): the Malgrange--Fay form will undergo a resulting modification but the main theorem that its exterior derivative is the pull-back of the Goldman structure should remain true. This implies that the difference of the Malgrange--Fay forms for the original bundle and the Hecke-modified  is the differential of a locally defined function. The question is then to explicitly characterize this function. 
\item[-] Fay's identified (formula (10) in \cite{Fay92}, proved in \cite{Fay92book}) the one form  with the holomorphic part of a suitable factorization of the analytic Ray-Singer torsion (see (10) in \cite{Fay92} and also Corollary to Theorem 1 in \cite{Quillen}). In view of a generalization to meromorphic connections, the tau function described above will provide the analog of a ``determinant'' for wild character varieties.

\item[-]  The Tyurin parametrization seems suited for studying the degeneration where the base curve develop one or more nodes. The reason is that the main tool, the Cauchy kernel, is written in terms of holomorphic differentials and third kind differentials, whose behaviour under degeneration is classically well studied. This may allow to extend results of \cite{DaskWent} who studied in detail the rank two case.
\end{itemize}

\appendix
\renewcommand{\theequation}{\Alph{section}.\arabic{equation}}

\section{Proof of Theorem \ref{thmmalgfay}}
\label{proofmain}

There are two main components in the definition \eqref{MalgFay}; the fundamental flat matrix $\Psi$ and the reference connection $\F_{\scr D}$. To compute the exterior derivative of $\Malg$ we need variational formulas for both.

\begin{lemma}
\label{lemmavarPsi}
The variational formula for $\Psi$ takes the form
\be
\delta \Psi({p})\Psi^{-1}({p})=\frac{1}{2\pi\i}\int_{q\in\Sigma}  \Psi_-({q})\delta J({q})J^{-1}({q})\Psi_-^{-1}({q})\Cauchy_\infty({q},{p})
.\label{eqvarPsi}
\ee
\end{lemma}
\noindent{\bf Proof.}
It follows from the jump condition $\Psi_+(p)=\Psi_-(p)J(p)$ and the normalization $\Psi(\infty)=\1$ that $\delta\Psi\Psi^{-1}$ satisfies the {\it additive} jump condition
\be
\delta\Psi_+\Psi^{-1}_+=\delta\Psi_-\Psi^{-1}_-+\Psi_-\delta J J^{-1}\Psi_-^{-1}, 
\ee
as well as the normalization $\delta\Psi(\infty)=0$. 
It follows from the Sokhotskii-Plemelji formula that the expression in the right side of \eqref{eqvarPsi} satisfies exactly the same conditions. Denote by $R(p)$ their difference; 
since the divisor $\scr D$ is not part of the variation, it is holomorphic at $\scr D$. 
Then it has poles only at the Tyurin divisor. At each point $q\in \scr T$, the product $R(p)\nf (p)$ is analytic because of Corollary~\ref{corollCauchy} and
Theorem~\ref{nomonothm} (part 1). Thus the rows of $R(p)$ are sections of $\E$ that vanish at $\infty$. Since $\E$ is not on the non-abelian theta divisor, it follows that $R(p)\equiv 0$. 
\QED

\begin{lemma}[Variational formul\ae] 
\label{lemmavarF}
The variation of $\F_{\scr D} ({p})$ is 
\be
\label{varF}
\delta\F_{\scr D}(p)=\sum_{t\in\scr T}\res{q=t}\Cauchy_\infty(p,q)\delta\Psi(q)\Psi^{-1}(q)\Cauchy_\infty(q,p)
\ee
and it is independent of the choice of $\scr D$.
\end{lemma}

\noindent {\bf Proof.}
Observe that, as per Theorem \ref{nomonothm} (part 1)  the variational formula depends only on the polynomial normal form $\nf$ and not the whole solution $\Psi$.
The independence from $\scr D$ is clear because $\F_{\scr D} = \F - \1\d \ln h_{\scr D}$ and the reference spinor $h_{\scr D}$ is not subject to deformations, $\delta h_{\scr D}=0$. 
Observe now that $\mathbf K(q,p):= \Psi^{-1}(q)\Cauchy_\infty (q,p) \Psi(p)$ has a simple pole at $q=p$, with residue identity and hence moduli-independent, and is holomorphic otherwise, with a zero at $p=\infty$.
Therefore $\delta \mathbf K(q,p)$ is an analytic function of $p\in\CC\setminus\Sigma$ with a zero at $\infty$, $\mathbb K(q,p=\infty)=0$, and on $\Sigma$ it satisfies (for fixed ${q}\in\CC\setminus\Sigma$) 
\be
\delta \mathbf K({q}, {p}_+) =\delta \mathbf  K({q}, {p}_-)  J({p})   + \mathbf K({q}, {p}_-) \delta J({p})\ ,\qquad
{p}\in\Sigma.
\ee
It follows again from Sokhotskii-Plemelji  formula,  along the same lines as  Lemma \ref{lemmavarPsi} that 
\be
\label{start}
\delta \mathbf K(q,p)  = \frac 1{2\pi\i}\int_{s\in\Sigma}\mathbf K(q,s_-) \delta J(s)J^{-1}(s) \mathbf K(s_-,p)
\ee
because  the expression on the right side satisfies the same jump condition and normalization at ${p}=\infty$.
Using the notation \eqref{jumpdef} and the Cauchy theorem the formula \eqref{start} can be written as a sum of residues at $p, q, \scr T$ as follows
\begin{align} 
\delta \mathbf K({q},{p}) 
&= \Psi^{-1}({q}) \Cauchy_\infty({q},{p}) \delta \Psi({p})- \Psi^{-1}({q})  \delta \Psi({q}) \Psi^{-1}({q}) \Cauchy_\infty({q},{p}) \Psi({p}) \nn\\
&\qquad+ \sum_{{t} \in \mathscr T} \res{s=t}
\Psi^{-1}({q})\Cauchy_\infty(q,s) \delta\Psi(s) \Psi^{-1} (s)  \Cauchy_\infty(s,p) \Psi(p)
\nn\\
&=\mathbf K(q,p) \Psi^{-1}(p) \delta\Psi(p) - \Psi^{-1}(q) \delta\Psi(q)\mathbf K(q,p) \nn\\
&\qquad+\sum_{{t} \in \mathscr T} \res{s={t}}
\Psi^{-1}(q)\Cauchy_\infty(q,s) \delta\Psi(s) \Psi^{-1} (s)  \Cauchy_\infty({s},p) \Psi(p).
\label{13}
\end{align}
Let $z= \zeta({p})$ and $w = \zeta({q})$ for a local coordinate $\zeta$; by a slight abuse of notation we shall indicate $\Psi(z)$ the evaluation of $\Psi$ at $p$, and similarly for $\mathbf K$ and for the point $q$. Taking the regular part at $w=z$ of the kernel $\mathbf K$ yields
\be
\lim_{w\to z} \le( \frac{\mathbf K(w,z)}{\d w} - \frac \1{w-z} \ri)= \Psi^{-1}(z)\F(z)\Psi(z)- \Psi^{-1}(z)\frac{\d}{\d z}\Psi(z).
\ee
We now take a variation of the moduli and obtain (we denote for brevity $'$ the derivative w.r.t $z$ in the chosen coordinate and observe that $\delta\mathbf K(w,z)$ is regular along $w=z$): 
\begin{align}
\delta \mathbf K(z,z) 
= \bigg[\Psi^{-1}(z) \F(z) \Psi(z), \Psi^{-1}(z) \delta \Psi(z)\bigg] - \delta( \Psi^{-1}(z)\Psi'(z)) + \Psi^{-1} (z)\delta \F(z) \Psi(z).
\label{24}
\end{align}
On the other hand, the limit $q\to p$ of \eqref{13}  gives (all terms below are evaluated at $p$ unless otherwise indicated)
\begin{align}
\delta\mathbf K(p,p)&=   \bigg[ \Psi^{-1} \F \Psi, \Psi^{-1} \delta \Psi \bigg]  - \Psi^{-1} (\delta \Psi)' + \Psi^{-1} \delta \Psi \Psi^{-1} \Psi'
+  \sum_{t \in \mathscr T} \res{s=t}
\Psi^{-1}({p})\Cauchy_\infty({p} ,s) \delta\Psi(s) \Psi^{-1} (s)  \Cauchy_\infty(s,p) \Psi(p)
\nn
\\
&=   \bigg[ \Psi^{-1} \F \Psi, \Psi^{-1} \delta \Psi \bigg]  - \delta(\Psi^{-1} \Psi')+ \Psi^{-1}(p)\le( \sum_{t \in \mathscr T} \res{s=t }
\Cauchy_\infty(p,s) \delta\Psi(s) \Psi^{-1} (s)  \Cauchy_\infty(s,p) \ri)\Psi(p).
\label{25}
\end{align}
Then \eqref{varF} follows by comparing  \eqref{24} and \eqref{25}.
\QED

\noindent {\bf  Proof of Theorem~\ref{thmmalgfay}.} 
\label{proofthmmalgfay}
We split the Malgrange--Fay differential as follows 
\be
\label{723}
2\pi\i\, \Malg = \int_\Sigma \tr \bigg(\Psi_-^{-1} \d \Psi_- \delta J J^{-1} \bigg) - \int_{\Sigma} \tr \bigg(\F_{\scr D} \Psi_-^{-1} \delta J J^{-1} \Psi_{_-}\bigg)  =: \Malg_0 - \Malg_1
\ee
We are going to compute separately $\delta \Malg_0$ and $\delta \Malg_1$ and we will show that
\begin{align}
\label{dtheta0}
\delta \Malg_0
&= -\frac 1 2 \Omega(\Sigma)   
+\frac 1 2\sum_{t\in \mathscr T}  \res{{p}=t} \tr \bigg(
  \d( \delta \Psi \Psi^{-1}) \wedge \delta \Psi \Psi^{-1}\bigg),
\\
\label{dPhi}
\delta \Malg_1 
&= \frac 12 \sum_{t\in \mathscr T} \res{p=t} \tr  \bigg(
\d\le(\delta \Psi\Psi^{-1}\ri)\wedge \delta \Psi \Psi^{-1} \bigg),
\end{align}
from which the theorem will follow.

\paragraph{Computation of $\delta \Malg_0$.}
For brevity we will use the notation $\Cartan{}:=\delta J J^{-1}$. We have
\begin{align}
\label{step1}
\delta \Malg_0=& 
 \int_\S \tr \bigg(
 -\Psi_-^{-1} \delta \Psi_- \Psi_-^{-1}  \d \Psi_- \wedge \Cartan{}
 + \Psi^{-1}_- \delta \d \Psi_- \wedge \Cartan{}  + \Psi_-^{-1} \d \Psi_- \delta\Cartan{} \bigg).
\end{align}
Now use the identity $\delta \Cartan{} = \Cartan{} \wedge \Cartan{} $ as well as the variational formul\ae
\begin{align}
 \delta \Psi({p}) &=\int_{q\in\S} \Psi_-({q}) \Cartan{(q)} \Psi_-^{-1}({q})\Cauchy_\infty({q},{p}) \Psi({p}),\nn
\\
\delta \d \Psi({p})
& = \delta \Psi({p}) \Psi({p})^{-1}\d \Psi({p})
+ \int_{q\in\S} \Psi_-({q}) \Cartan{(q)} \Psi_-^{-1}({q}) \d_{p}\Cauchy_\infty({q},{p}) \Psi({p}),
\label{deltaPsi}
\end{align}
the first of which is the content of Lemma \ref{lemmavarPsi} and the second is a consequence. Inserting \eqref{deltaPsi} into \eqref{step1} we obtain 
\be
\label{btt0}
\delta \Malg_0=\bs\chi + \int_{\Sigma} \tr \bigg(\Psi_-^{-1} \d  \Psi_- \Cartan{} \wedge \Cartan{}\bigg) 
\ee
where we denote
\be
\bs\chi:=\int_{{p} \in \Sigma} \int_{{q} \in \Sigma} \tr \bigg(
    \Psi_-({q} ) \Cartan{(q)}\Psi_-({q} )^{-1} \d_{p} \Cauchy_\infty({q} ,{p} _-) \wedge\Psi_-({p} ) \Cartan{(p)}\Psi({p} )^{-1}_-\bigg).
\ee
Here the subscript ${p}_-$ in the Cauchy kernel indicates that the integration is taken on the boundary value from the right of $\S$, which is important to specify since the integral is  a singular integral.

To handle the term $\bs\chi$ in \eqref{btt0} it is convenient to introduce the {\it fundamental bi-differential} $\Berg(q,p)$ on $\CC$, which is the unique symmetric scalar bi-differential with a double pole at ${p}={q}$ with unit bi-residue and trivial $\alpha$-periods \cite[Ch. III]{Fay73}.
The only property we are going to use is its symmetry $\Berg(q,p)=\Berg(p,q)$. Then, we note that $\d_{p}\Cauchy_\infty({q},{p}) = \Berg ({q},{p})\1 + \mathbf H({q},{p})$, where $\mathbf H(q,p)$ is a  matrix bi-differential regular along ${q}={p}$. Thus $\bs\chi=\bs\chi_1+\bs\chi_2$ with
\begin{align}
\bs\chi_1:={}& \int_{{p} \in \Sigma} \int_{{q} \in \Sigma} \tr \bigg(
    \Psi_-({q} ) \Cartan{(q)} \Psi_-({q} )^{-1}\Berg({q} ,{p} _-) \wedge\Psi_-({p} ) \Cartan{(p)}\Psi({p} )^{-1}_-\bigg), \nn\\
\bs\chi_2:={}& \int_{{p} \in \Sigma} \int_{{q} \in \Sigma} \tr \bigg(
    \Psi_-({q} ) \Cartan{(q)}\Psi_-({q} )^{-1} \mathbf H({q} ,{p} ) \wedge\Psi_-({p} ) \Cartan{(p)}\Psi({p} )^{-1}_-\bigg).
    \label{defA}
\end{align}

\begin{proposition}
\label{propa4}
The singular integral $\bs\chi_1$ can be evaluated as
\be
\label{chi1}
\bs\chi_1= - \int_{\Sigma} \tr \bigg( \Psi_-^{-1} \d  \Psi_-\Cartan{}  \wedge \Cartan{}\bigg) -\frac 1 2 \Omega(\Sigma).
\ee
\end{proposition}

We only sketch the proof, because the computation can be handled following identical steps (which we now review) as the proof of Theorem~2.1 in~\cite{Bertocorr}.
Indeed, we observe that the only properties used in the proof in loc. cit. are that the kernel of the integration has a singularity of type $\d z\d w(z-w)^{-2}$ on the diagonal and it is symmetric in the two variables. Both properties are satisfied by the fundamental bidifferential $\Berg$ (Ch. III of \cite{Fay73}).

We start by introducing a small disk $\mathbb D_\epsilon$ on $\CC$, centered at $\infty$ and of radius $\epsilon$ in some coordinate chart, and separate the integration over $\S\cap\DD_\epsilon$ and $\S\setminus\DD_\epsilon$.
Namely, we decompose $\bs\chi_1=A_\epsilon+B_\epsilon+C_\epsilon$, where
\begin{align}
A_\epsilon:={}&\int_{{p} \in \Sigma\setminus\DD_\epsilon} \int_{{q} \in \Sigma\setminus\DD_\epsilon} \tr \bigg(
    \Psi_-({q} ) \Cartan{(q)} \Psi_-({q} )^{-1}\Berg({q} ,{p} _-) \wedge\Psi_-({p} ) \Cartan{(p)}\Psi({p} )^{-1}_-\bigg),
\\
B_\epsilon:={}&\int_{{p} \in \S\cap\DD_\epsilon} \int_{{q} \in \Sigma\setminus\DD_\epsilon} \tr \bigg(
    \Psi_-({q} ) \Cartan{(q)} \Psi_-({q} )^{-1}\Berg({q} ,{p} _-) \wedge\Psi_-({p} ) \Cartan{(p)}\Psi({p} )^{-1}_-\bigg)
   \nn \\
    &\quad
    +    \int_{{p} \in \Sigma\setminus\DD_\epsilon} \int_{{q} \in \S\cap\DD_\epsilon} \tr \bigg(
    \Psi_-({q} ) \Cartan{(q)} \Psi_-({q} )^{-1}\Berg({q} ,{p} _-) \wedge\Psi_-({p} ) \Cartan{(p)}\Psi({p} )^{-1}_-\bigg),
\\
C_\epsilon:={}&\int_{{p} \in \Sigma\cap\DD_\epsilon} \int_{{q} \in \Sigma\cap\DD_\epsilon} \tr \bigg(
    \Psi_-({q} ) \Cartan{(q)} \Psi_-({q} )^{-1}\Berg({q} ,{p} _-) \wedge\Psi_-({p} ) \Cartan{(p)}\Psi({p} )^{-1}_-\bigg).
\end{align}

We compute separately these three contributions to $\bs\chi_1$; each of them is treated in the same way as in loc. cit.
Since $\S\setminus\DD_\epsilon$ is a disjoint union of smooth arcs without self-intersections, for the computation of $A_\epsilon$ we can rely on the following result, which follows from the classical theory of singular integrals (see e.g. \cite{Gakhov}).

\begin{lemma}
\label{lemmasame}
Let $\gamma$ be an oriented smooth arc in $\CC$ without self-intersection and let $\varphi: \gamma \times\gamma \to \C$ be a function which is locally analytic in each variable and such that $\varphi(z,w)= - \varphi(w,z)$.  
Then 
\be
\int_{z\in\gamma} \int_{w\in \gamma} \Berg(z_-,w) \varphi(w,z)  = -\frac 1 2 \int_\gamma \d_w \varphi(w,z)\bigg|_{w=z}.
\ee
\end{lemma}
This gives
\begin{align}
A_\epsilon&=- \frac 1 2 \int_{\Sigma\setminus \DD_\epsilon} \tr \bigg( \d_{q} \le( \Psi_-({q} ) \Cartan{(q)}\Psi_-({q} )^{-1} \ri) \wedge\Psi_-({p} ) \Cartan{(p)}\Psi({p} )^{-1}_-\bigg) \bigg|_{{q} ={p} }\nn
\\
&=-\frac 1 2 \int_{\Sigma\setminus \S_\epsilon} \tr \bigg( \d  \Psi_-\Cartan{}  \wedge \Cartan{}\Psi^{-1}_-  -  \Cartan{}\Psi_-^{-1}\d  \Psi_-  \wedge \Cartan{}  + \d\,  \Cartan{} \wedge \Cartan{}  \bigg)
\nn\\
& = -\int_{\Sigma\setminus\S_\epsilon} \tr \bigg(\Psi_-^{-1} \d  \Psi_-\Cartan{}  \wedge \Cartan{}   + \d\,  \Cartan{} \wedge \Cartan{}  \bigg).
\end{align}
Thus $A_\epsilon$ clearly admits the limit
\be
\lim_{\epsilon\to 0_+}A_\epsilon=-\int_{\Sigma} \tr \bigg( \Psi_-^{-1} \d  \Psi_-\Cartan{}  \wedge \Cartan{}   + \d\,  \Cartan{} \wedge \Cartan{}  \bigg).
\ee
Finally, the same reasoning in \cite{Bertocorr} applies to show that $B_\epsilon=0$ because of skew-symmetry and that $C_\epsilon$ yields the localized contribution
\be
\lim_{\epsilon\to 0_+}C_\epsilon=-\frac 1 2\Omega(\Sigma).
\ee
This part of the computation is identical to the one in loc. cit. because it can be carried out entirely in local coordinate and requires only the symmetry of $\Berg$ together with its singular behaviour along the diagonal.
Finally, since $J({p})$ in \eqref{jumpdef} is constant on each arc of $\Sigma$, we have $\d\Cartan{}=0$ and \eqref{chi1} follows.

To complete the computation of 
\be
\delta\Malg_0=-\frac 12\Omega(\S)+\bs\chi_2,
\ee
cf. \eqref{btt0}, it remains to compute $\bs\chi_2$.
Note that from \eqref{jumppsi} we obtain $\Psi_- \Cartan{} \Psi^{-1}_- =\Psi_- \delta J J^{-1} \Psi^{-1}_-= \Delta_\Sigma(\delta \Psi \Psi^{-1})$, where $\Delta_\Sigma$ is the jump operator across $\Sigma$, see \eqref{jumpdef}; hence we can use Cauchy residue theorem to write 
\begin{align}
\bs\chi_2 &= \int_{p\in \Sigma} \int_{{q} \in \Sigma} \tr \bigg(
   \Delta_\Sigma (\delta \Psi(q) \Psi^{-1}({q} ))  \mathbf H({q} ,p) \wedge \Delta_\Sigma (\delta \Psi (p)\Psi^{-1}(p))\bigg)
  \nn \\
   &= \sum_{t\in \mathscr T} \res{{q} =t}   \int_{p\in \Sigma}  \tr \bigg(
  \delta \Psi(q) \Psi^{-1}({q} )  \mathbf H({q} ,p) \wedge \Delta_\Sigma (\delta \Psi(p) \Psi^{-1}(p))\bigg)
  \label{231}
\end{align}
To proceed, we now show that
\be
\label{HtoB}
\res{{q} =t} \delta \Psi({q} ) \Psi^{-1}({q} )  \mathbf H({q} ,{p}) = - \res{{q} =t} \delta \Psi({q} ) \Psi^{-1}({q} ) \Berg({q} ,{p}).
\ee
To this end we note that, by construction of the Cauchy kernel, $\Psi^{-1}({q}) \Cauchy_\infty({q},{p})$ is a differential with respect to ${q}$ {\it without poles} at the Tyurin divisor $\scr T$ (this follows from Theorem-Definition~\ref{defCauchy} and Theorem~\ref{nomonothm}). Therefore
\begin{align}
\forall t\in \scr T:\ \ \ \ \res{{q} =t} \delta \Psi({q} ) \Psi^{-1}({q} ) \Cauchy_\infty({q} ,{p})=0 \ \ \ \Rightarrow \ \ \ \res{{q} =t} \delta \Psi({q} ) \Psi^{-1}({q} ) \d_{p}\Cauchy_\infty({q} ,{p})=0.
\end{align}
Recalling that $\d_{p}\Cauchy_\infty({q} ,{p}) = \Berg(q,p)\1 + \mathbf H(q,p)$ we obtain \eqref{HtoB}. 
Therefore  in \eqref{231} we can replace the matrix kernel $\mathbf H$ with the (scalar) kernel $-\Berg$, use Cauchy theorem again  and obtain 
\begin{align}
\bs\chi_2  =
- \sum_{{t}\in \mathscr T} \res{p={t}}  \sum_{{r}\in \mathscr T} \res{{q} ={r}}  \tr \bigg(
  \delta \Psi(q) \Psi^{-1}({q} )\wedge \delta \Psi(p) \Psi^{-1}(p)\bigg) \Berg ({q} ,z) .
\end{align}
In the double sum over the points of the Tyurin divisor $\scr T$  only the diagonal part contributes because of the skew-symmetry of the expression.  
To compute the iterated residue at $t\in\scr T$, we work in a local coordinate $\zeta$ centered at $t$ and denote 
\be
\varphi(w,z):= \tr \bigg(
  \delta \Psi(w) \Psi^{-1}({w} )\wedge \delta \Psi(z) \Psi^{-1}(z)\bigg),
\ee
with $w=\zeta(q)$, $z=\zeta(p)$. This expression is jointly analytic in the punctured disks around $z=0$ and $w=0$, and skew symmetric, $\varphi(z,w)=-\varphi(w,z)$.
Using Lemma \ref{lemmasame}
we obtain 
\be
\res{z=0} \res{w=0} \Berg(z,w) \varphi(z,w)\d z \d w = -\frac 1 2 \res{z=0} \d_w \varphi(z,w)\big|_{w=z}.
\ee
From this we finally obtain 
\begin{align}
\bs\chi_2
&  = 
- \sum_{t \in \mathscr T} \res{{p}=t }  \res{z =t }  \tr \bigg(
  \delta \Psi \Psi^{-1}({p} ) \Berg ({p} ,z) \wedge \delta \Psi \Psi^{-1}(z)\bigg)   =\frac 1 2\sum_{t \in \mathscr T}  \res{{p}=t } \tr \bigg(
  \d_{p}( \delta \Psi \Psi^{-1}) \wedge \delta \Psi \Psi^{-1}\bigg),
\end{align} 
and equation \eqref{dtheta0} follows.

\paragraph{Computation of $\delta \Malg_1$.}
Using again that $\Psi_-^{-1} \Cartan{} \Psi_{_-} = \Delta_\Sigma ( \delta \Psi \Psi^{-1})$ and Cauchy residue theorem, we obtain
\begin{align}
\Malg_1& = \int \tr \bigg( \F_{\scr D}  \Psi_- \Cartan{}\Psi_-^{-1} \bigg) =\sum_{t\in \mathscr T} \res{{p}=t} \tr  \bigg( \F_{\scr D} \delta \Psi \Psi^{-1}\bigg).
\end{align}
In the use of Cauchy theorem we have noted that $\F_{\scr D} ({p})$ has a simple pole at ${p}=\infty$ but $\delta \Psi \Psi^{-1}$ has a simple zero such that there is no contribution to the residue theorem at ${p}=\infty$. 
We remark, in passing, that this term $\Malg_1$ does not have a counterpart in the genus zero setting of \cite{BertoMalg,Bertocorr}.

We can use the variational formula for $\F_{\scr D}$ (see Lemma \ref{lemmavarF}) to complete the computation; 
\begin{align}
\label{t1}
\delta \Malg_1  &= \sum_{t \in \mathscr T} \res{{p}=t} \le[\tr  \bigg( \delta \F_{\scr D}  \wedge \delta \Psi \Psi^{-1}\bigg) +  \tr  \bigg( \F_{\scr D}   \delta \Psi \Psi^{-1} \wedge \delta \Psi \Psi^{-1}\bigg) \ri]\nn\\
&\mathop{=}^{\eqref{varF}} 
\sum_{t \in \mathscr T} \res{{p}=t} \sum_{s \in \mathscr T} \res{{r} =s }\tr  \bigg( 
 \Cauchy_\infty({p} ,{r} ) \delta\Psi({r} ) \Psi^{-1} ({r} )  \Cauchy_\infty({r} ,{p}) \wedge \delta \Psi({p} ) \Psi({p} )^{-1}\bigg)
 +
 \sum_{ t \in \mathscr T} \res{{p} =t } \tr  \bigg( \F_{\scr D}  \delta \Psi \Psi^{-1} \wedge \delta \Psi \Psi^{-1}\bigg).
\end{align}
Denote by $(D)$ the double sum in \eqref{t1}. We observe (using the cyclic property of the trace) that the argument is odd under the exchange $p\leftrightarrow r$ with a simple pole at $p=r$. Then, using the same logic used in the previous part of the proof we easily obtain
\begin{align}
\nn
(D)&=-\frac 12 \sum_{t\in \mathscr T}  \res{{p}=t}\res{{r}=p}  
\tr  \bigg( 
 \Cauchy_\infty({p} ,{r} ) \delta\Psi({r} ) \Psi^{-1} ({r} )  \Cauchy_\infty({r} ,{p} ) \wedge \delta \Psi({p} ) \Psi^{-1}({p} )\bigg)
 \\
&= -\frac 12 \sum_{t \in \mathscr T} \res{{p} =t } \tr  \bigg(-
\d \le(\delta \Psi\Psi^{-1}\ri)\wedge \delta \Psi \Psi^{-1} 
+\F_{\scr D} \delta \Psi\Psi^{-1}\wedge \delta \Psi\Psi^{-1} - \delta \Psi\Psi^{-1} \F_{\scr D} \wedge \delta \Psi\Psi^{-1}
\bigg).
\label{118}
\end{align}
Note that in the above formula $\tr (\F_{\scr D} \delta \Psi \Psi^{-1}\wedge\delta \Psi\Psi^{-1})$ is actually a  well defined (local) differential with respect to change of coordinates because $\tr (\delta \Psi \Psi^{-1} \wedge \delta \Psi\Psi^{-1})=0$.

Then
\begin{align}
(D) =\frac 12 \sum_{t\in \mathscr T} \res{t} \tr  \bigg(
\d\le(\delta \Psi\Psi^{-1}\ri)\wedge \delta \Psi \Psi^{-1} \bigg) 
   - \sum_{t\in \mathscr T} \res{t} \tr \bigg(\F_{\scr D} \delta \Psi\Psi^{-1}\wedge \delta \Psi\Psi^{-1} \bigg)
   \label{diamante}
\end{align}
Finally,  plugging \eqref{diamante} in \eqref{t1} we obtain \eqref{dPhi}. 
Plugging \eqref{dtheta0} and \eqref{dPhi} in the differential of the relation \eqref{723},   $2\pi\i \delta \Malg = \delta \Malg_0 - \delta \Malg_1$, gives the theorem.\QED

\paragraph*{Data availability statement.}
Data sharing not applicable to this article as no datasets were generated or analysed during the current study.

\end{document}